\documentclass[11pt]{article}
\usepackage[a4paper,margin=1.12in]{geometry}
\usepackage{amsmath,amssymb,amsthm,mathtools}
\usepackage{mathrsfs}
\numberwithin{equation}{section}
\usepackage[T1]{fontenc}
\usepackage{lmodern}
\usepackage{microtype}
\usepackage{xcolor}
\usepackage{tikz}
\definecolor{linkblue}{rgb}{0.05,0.15,0.5}
\usepackage[colorlinks=true,linkcolor=linkblue,citecolor=green!60!black,urlcolor=linkblue,pdfborder={0 0 0}]{hyperref}
\hypersetup{bookmarksnumbered=true,bookmarksopen=true}
\newtheorem{theorem}{Theorem}[section]
\newtheorem{lemma}[theorem]{Lemma}
\newtheorem{proposition}[theorem]{Proposition}
\newtheorem{corollary}[theorem]{Corollary}
\theoremstyle{definition}
\newtheorem{definition}[theorem]{Definition}
\newtheorem{example}[theorem]{Example}
\newtheorem{remark}[theorem]{Remark}

\newlength{\condwd}
\newlength{\condind}
\newcommand{\condition}[2]{\par\noindent\settowidth{\condwd}{\textbf{#1}}\ifdim\condwd<1.9em \setlength{\condwd}{1.9em}\fi\setlength{\condind}{\dimexpr1.2em+\condwd+0.6em\relax}\hangindent=\condind\hangafter=1\hspace*{1.2em}\makebox[\condwd][l]{\textbf{#1}}\hspace*{0.6em}#2\par\hangindent=0pt\hangafter=1\vspace{0.2em}}
\newcommand{\subcondition}[2]{\par\noindent\settowidth{\condwd}{\textbf{#1}}\ifdim\condwd<0.9em \setlength{\condwd}{0.9em}\fi\setlength{\condind}{\dimexpr2.7em+\condwd+0.6em\relax}\hangindent=\condind\hangafter=1\hspace*{2.7em}\makebox[\condwd][l]{\textbf{#1}}\hspace*{0.6em}#2\par\hangindent=0pt\hangafter=1\vspace{0.1em}}
\newcommand{\cl}[1]{\langle #1\rangle}
\newcommand{\G}{\mathcal G}
\newcommand{\K}{\mathcal{K}}
\newcommand{\fin}{\subseteq_{\mathrm{fin}}}
\newcommand{\ua}{\mathord{\uparrow}}
\newcommand{\da}{\mathord{\downarrow}}
\newcommand{\dda}{\mathord{\mbox{\makebox[0pt][l]{\raisebox{-.4ex}{$\downarrow$}}$\downarrow$}}}
\newcommand{\dua}{\mathord{\mbox{\makebox[0pt][l]{\raisebox{.4ex}{$\ua$}}$\ua$}}}
\newcommand{\set}[2]{\{#1\mid#2\}}
\newcommand{\UB}{\operatorname{UB}}
\newcommand{\mub}{\operatorname{mub}}
\newcommand{\Fun}{\mathscr F}
\newcommand{\Funs}{\mathscr F_{\mathrm s}}
\begin{document}
\begin{center}
{\large\bfseries Lawson compact domains and generalized closure spaces\par}
\vspace{1.5em}
{\normalsize Longchun Wang\par}
\vspace{0.5em}
{\small School of Mathematical Sciences, Qufu Normal University,\\
Qufu, Shandong 273165, China.\\
Corresponding author email: longchunw@163.com.}
\end{center}
\begin{quote}
\noindent\textbf{Abstract.}
We study Lawson compact domains and certain subclasses within
the framework of interpolative generalized closure spaces
(IG-closure spaces). We give two conditions on an IG-closure
space that characterize Lawson compactness of its domain of
generalized directed sets, yielding a representation of Lawson
compact domains. We then characterize domains with property $M$, as
well as domains admitting a basis with property $M$. We also
represent Lawson compact algebraic domains by dense IG-closure
spaces satisfying two additional conditions.
Under local consistency, these conditions can be simplified for
Lawson compact algebraic $L$-domains. For these classes of
domains, we establish corresponding categorical equivalences.
In addition, we construct a pointed Lawson compact $L$-domain in
which every basis contains a pair with infinitely many minimal
upper bounds in that basis. This answers a question of Jiang, Li
and Zhou in the negative.
\end{quote}

\noindent\textbf{Mathematics Subject Classification (2020):} 06B35, 06A15, 54A05, 18B35.\\
\textbf{Keywords:} Domain theory, Lawson compact domain, interpolative
generalized closure space, 
algebraic $L$-domain, categorical equivalence.
\section{Introduction}
The Lawson topology on a domain is generated by its Scott and
lower topologies. Lawson compactness has been studied in connection
with coherence and function spaces \cite{kouluo,liangkeimel}.
Several characterizations of this property are known. For continuous dcpos, Jung characterized
Lawson compactness by Scott compactness of the whole space and of
intersections of any two Scott compact upper sets
\cite[Lemma~4.17]{jung}. For well-filtered dcpos, Jia, Jung and Li
characterized Lawson compactness by Scott compactness of the whole
space and of intersections of any two principal upper sets
\cite[Theorem~3.4]{jiajungli}. For further studies on Lawson compactness, see
\cite{goubaultjung,lawsonxi,luanli}.

Lawson compactness can also be characterized in terms of a basis.
Let $D$ be a domain with basis $\mathcal B$. By
\cite[Corollary~III-5.13]{gierz}, $D$ is Lawson compact if and
only if it is finitely generated and, for any
$a_1,a_2,b_1,b_2\in\mathcal B$ with $b_i\ll a_i$ for
$i=1,2$, there exists a finite set $F\subseteq\mathcal B$ such
that
\[
\ua a_1\cap\ua a_2
\subseteq\ua F
\subseteq\ua b_1\cap\ua b_2.
\]
The condition expressed by this formula is called property
$M$ with respect to $\mathcal B$ in \cite{gierz}.
The term property $M$ is also used for a property of posets
in \cite[Section~1.4]{jung}. A poset $P$ has property $M$ in
this sense if it has property $m$ and every finite subset of $P$
has finitely many minimal upper bounds in $P$. In this paper,
we use property $M$ in this latter sense.

A domain is Lawson compact if it admits a basis with
property $M$ \cite[Remark~2.8]{jiang}. For an algebraic domain,
this condition is also necessary, as its compact elements form a
basis with property $M$ whenever the domain is Lawson compact
\cite[Corollary~4.19]{jung}; see also \cite{abramskyjung,gierz}.
Jiang, Li and Zhou gave a logical representation of $L$-domains
admitting a basis with property $M$ \cite[Section~4.1]{jiang}.
They also asked whether every Lawson compact $L$-domain admits
such a basis \cite[Section~5]{jiang}. Their question can be stated
as follows.
\condition{(Q1)}{Does every Lawson compact $L$-domain admit a
basis with property $M$?}

 In this paper, we aim to study Lawson compact domains using
interpolative generalized closure spaces (IG-closure spaces). Representations based on closure
spaces have a long history. A classical example is that every
algebraic lattice is isomorphic to the lattice of closed subsets
of a set equipped with an algebraic closure operator
\cite[Theorem~5.8]{burris}.  Such representations have been developed for algebraic domains
and algebraic $L$-domains \cite{guoli,wuguoli}, continuous
lattices \cite{liwangyao}, FS-domains \cite{yaoli-fs}, and
stably continuous semilattices \cite{yao}. Lawson compact
$L$-domains have also been represented by information systems
\cite{jianginfo}. In the algebraic case, logical and topological
representations have been obtained in \cite{wangli2020} and
\cite{wangzhouli-top}, respectively.

Wang and Li \cite{wangli} represented domains and several of their
subclasses by IG-closure spaces. A quantale-valued extension was studied in
\cite{wuyaoli}. For an IG-closure space $X$, its generalized
directed sets, ordered by inclusion, form a domain $\G(X)$ with
canonical basis $\mathcal B_X=\set{\cl x}{x\in X}$. Every domain is, up to order
isomorphism, of this form. This raises the following question.
\condition{(Q2)}{How can Lawson compact domains be characterized
in terms of IG-closure spaces?}
For a subset $S\subseteq X$, put
$\mathcal B_S=\set{\cl s}{s\in S}$. We
also consider the following questions.

\condition{(Q3)}{Under what conditions on an IG-closure space $X$
and a subset $S\subseteq X$ is $\mathcal B_S$ a basis of $\G(X)$
with property $M$?}

\condition{(Q4)}{Can every domain admitting a basis with
property $M$ be represented as $\G(X)$ with $\mathcal B_S$
as such a basis?}

After recalling the necessary preliminaries in Section~\ref{sec:pre},
we answer \textup{(Q2)} in Section~\ref{sec:lc}.
Theorem~\ref{thm:main} characterizes Lawson compactness of $\G(X)$
by two conditions, \textup{(C1)} and \textup{(C2)}, which are independent
(Examples~\ref{ex:C1notC2} and~\ref{ex:C2notC1}). This yields a
representation of Lawson compact domains by IG-closure spaces
satisfying these conditions (Theorem~\ref{thm:repr_domain}).

In Section~\ref{sec:m}, we answer \textup{(Q3)}. We introduce
a condition \textup{(CM)} on an IG-closure space $X$ and show
that it is equivalent to the canonical basis $\mathcal B_X$
having property $M$ (the case $S=X$). For an arbitrary subset
$S\subseteq X$, we then characterize when $\mathcal B_S$ is a
basis of $\G(X)$ with property $M$, namely by conditions
\textup{(B)} and \textup{(M$_S$)} (Lemmas~\ref{lem:subbasis}
and~\ref{lem:ms}). Combined with Theorem~\ref{thm:repr_domain},
this yields an affirmative answer to \textup{(Q4)} (Theorem~\ref{thm:selected}).

These results allow us to compare Lawson compactness of
$\G(X)$ with property $M$ of its bases. Condition \textup{(CM)}
implies \textup{(C1)} and \textup{(C2)}, but the converse fails.
In Example~\ref{ex:CMfails}, $\G(X)$ has a basis with property
$M$, although its canonical basis does not. In
Example~\ref{ex:nobasis}, $\G(X)$ is Lawson compact but has no
basis with property $M$. The domain constructed in the latter
example is a pointed $L$-domain, giving a negative answer to
\textup{(Q1)}.

In Section~\ref{sec:alg}, we consider algebraic domains. For a dense
IG-closure space $X$, we characterize Lawson compactness of $\G(X)$
by conditions \textup{(A1)} and \textup{(A2)}, which involve only
elements $x\in X$ satisfying $x\in\cl x$
(Theorem~\ref{thm:algebraic}). This yields a representation of
Lawson compact algebraic domains by dense IG-closure spaces
satisfying these conditions (Theorem~\ref{thm:alg_domain}). If $X$
is also locally consistent, Lawson compactness of $\G(X)$ admits a
characterization in terms of the minimal elements and minimal upper
bounds of pairs in $\K(\G(X))$ (Theorem~\ref{prop:algL-mub}). This
gives a representation of Lawson compact algebraic $L$-domains
(Theorem~\ref{thm:algL-mub-repr}). Although the canonical IG-closure space of a Lawson compact
algebraic domain need not satisfy \textup{(CM)}
(Example~\ref{ex:CMfails}), every such domain has a dense
representing IG-closure space satisfying \textup{(CM)}, obtained
from its compact elements (Theorem~\ref{thm:algCM}).

Finally, in Section~\ref{sec:cat}, we restrict the categorical
equivalence of Wang and Li \cite{wangli} to obtain equivalences
for Lawson compact domains, Lawson compact algebraic domains,
Lawson compact algebraic $L$-domains, and domains that themselves
have property $M$. We also establish a categorical equivalence between
IG-closure spaces equipped with suitable subsets and domains
equipped with bases having property $M$.

\section{Preliminaries}\label{sec:pre}

\subsection{Domains and the Lawson topology}

We first recall the standard notions of domain theory used below; see \cite{gierz} for details.

Let $P$ be a poset. For $A\subseteq P$, define
\[
\ua A=\set{p\in P}{a\le p\text{ for some }a\in A},
\qquad
\da A=\set{p\in P}{p\le a\text{ for some }a\in A}.
\]
For $x\in P$, write $\ua x=\ua\{x\}$ and $\da x=\da\{x\}$. When the ambient poset needs to be specified, we write $\ua_P A$ and $\da_P A$.
A subset $U$ of a poset $P$ is \emph{finitely generated} if $U=\ua F$ for $F\fin P$, and an \emph{upper set} if $U=\ua U$. A subset $V$ of $P$ is a \emph{lower set} if $V=\da V$.
A subset $E$ of $P$ is \emph{directed} if $E\ne\varnothing$ and, for every $x,y\in E$, there exists $z\in E$ such that $x\le z$ and $y\le z$. A poset is a \emph{dcpo} if every directed subset has a supremum. 

Let $D$ be a dcpo. For a subset $E\subseteq D$ that has a supremum, we write $\bigvee E$ for that supremum. For $x,y\in D$, we write $x\ll y$ if, for every directed subset $E\subseteq D$, $y\le\bigvee E$ implies $x\le e$ for some $e\in E$. For $x\in D$, put
\[
\dda x=\set{y\in D}{y\ll x},
\qquad
\dua x=\set{y\in D}{x\ll y}.
\]
A subset $\mathcal B$ of a dcpo $D$ is a \emph{basis} if $\mathcal B\cap\dda x$ is directed and $\bigvee(\mathcal B\cap\dda x)=x$ for every $x\in D$. A \emph{domain} is a dcpo with a basis. We write $F\fin A$ if $F$ is a finite, possibly empty, subset of $A$, and $F\ll x$ if $w\ll x$ for every $w\in F$.

 Let $D$ be a domain with basis $\mathcal B$. For every $x\in D$ and every $F\fin D$ with $F\ll x$ there exists $b\in \mathcal B$ such that $F\ll b\ll x$; see \cite[Theorem I-1.9]{gierz}. We refer to this as \emph{basis interpolation}. An element $k$ is \emph{compact} if $k\ll k$; $\K(D)$ denotes the set of all compact elements. Every basis contains $\K(D)$. A domain is \emph{algebraic} if $\K(D)$ is a basis. An \emph{$L$-domain} is a domain in which every principal ideal $\da x$ is a complete lattice under the induced order. A domain is \emph{pointed} if it has a least element.

\begin{definition}
Let $D$ and $D'$ be dcpos. For a map $f\colon D\to D'$ and $A\subseteq D$, write $f(A)=\set{f(a)}{a\in A}$. The map $f$ is \emph{Scott continuous} if it is monotone and
\[
f\Bigl(\bigvee E\Bigr)=\bigvee f(E)
\]
for every directed subset $E$ of $D$.
\end{definition}
\begin{definition}
Let $D$ be a dcpo. The \emph{Scott topology} $\sigma(D)$ consists of all upper sets $V\subseteq D$ for which $\bigvee E\in V$ implies $E\cap V\ne\varnothing$ for every directed $E\subseteq D$. The \emph{lower topology} $\omega(D)$ is generated by the subbasis $\set{D\setminus\ua x}{x\in D}$. The \emph{Lawson topology} $\lambda(D)$ is the smallest topology on $D$ containing both $\sigma(D)$ and $\omega(D)$. A dcpo $D$ is \emph{Lawson compact} if $(D,\lambda(D))$ is a compact topological space.
\end{definition}

\begin{lemma}[{\cite[Corollary III-5.13]{gierz}}]\label{lem:lawson}
Let $D$ be a domain with basis $\mathcal B$. Then $D$ is Lawson compact if and only if:
\begin{enumerate}
\item there exists $H\fin D$ with $D=\ua H$;
\item for any $a_1,a_2,b_1,b_2\in \mathcal B$ with $b_i\ll a_i$ $(i=1,2)$, there exists $F\fin \mathcal B$ such that
\[
\ua a_1\cap\ua a_2\subseteq\ua F\subseteq\ua b_1\cap\ua b_2.
\]
\end{enumerate}
\end{lemma}
The second condition of the above lemma is called property $M$ with
respect to $\mathcal B$ in \cite{gierz}. In this paper, the term
property $M$ is used in a different sense (Definition~\ref{def:property-m}).

\subsection{Generalized closure spaces}
For a set $X$, let $\mathcal P(X)$ denote its power set. We recall from \cite{wangli} the notions of generalized closure spaces used below.
\begin{definition}
A \emph{generalized closure space} is a set $X$ together with a monotone operation $\cl{\cdot}:\mathcal P(X)\to\mathcal P(X)$ satisfying $\cl{\cl A}\subseteq\cl A$ for every $A\subseteq X$. We write $\cl x$ for $\cl{\{x\}}$. An \emph{interpolative generalized closure space}, or \emph{IG-closure space}, is a generalized closure space satisfying the following condition:
\condition{(In)}{Whenever $x\in X$ and $M\fin\cl x$, there exists $y\in\cl x$ such that $M\subseteq\cl y$.}
\end{definition}

Note that $A\subseteq\cl B$ implies $\cl A\subseteq\cl B$. In particular, $x\in\cl y$ implies $\cl x\subseteq\cl y$. Taking $M=\varnothing$ in \textup{(In)} shows that $\cl x\neq\varnothing$ for every $x\in X$.
\begin{definition}\label{def:gd}
A subset $U$ of an IG-closure space $X$ is \emph{generalized directed} if, for every $M\fin U$, there exists $x\in U$ with $M\subseteq\cl x\subseteq U$. Write $\G(X)$ for the family of all generalized directed subsets of $X$, ordered by inclusion.
\end{definition}
\begin{lemma}\label{lem:gd}
Let $X$ be an IG-closure space and let $U\in\G(X)$. Then
\begin{enumerate}
\item $U\ne\varnothing$;
\item $\cl x\subseteq U$ for every $x\in U$.
\end{enumerate}
\end{lemma}
\begin{proof}
(1) Since $\varnothing\fin U$ and $U$ is generalized directed, there exists $x\in U$. Thus $U\ne\varnothing$.

(2) Let $x\in U$. Since $\{x\}\fin U$ and $U$ is generalized directed, there exists $y\in U$ with $x\in\cl y\subseteq U$. Hence $\cl x\subseteq\cl y\subseteq U$, because $a\in\cl b$ implies $\cl a\subseteq\cl b$.
\end{proof}
\begin{definition}[{\cite[Definitions~4.1--4.3]{wangli}}]\label{def:local}
Let $X$ be a generalized closure space. For $x,y\in X$, let $x\preceq y$ mean $\cl x\subseteq\cl y$, and for $A\subseteq X$ put
\[
\gamma(A)=\set{x\in X}{x\preceq a\text{ for some }a\in A}.
\]
We call $X$ \emph{pre-ordered} if the operation $\cl\cdot$ is compatible with $\preceq$, in the sense that
\begin{equation}\label{eq:compat}
\cl A=\cl{\gamma(A)}\subseteq\gamma(A)
\qquad\text{for every }A\subseteq X.
\end{equation}
A pre-ordered IG-closure space $X$ is said to be \emph{locally consistent} if, whenever $x\in X$ and $M\fin\cl x$, there exists $y\in\cl x$ such that $M\subseteq\gamma(\{y\})\subseteq\cl x$ and $\gamma(\{y\})\subseteq\gamma(\{z\})$ for every $z\in X$ satisfying $M\subseteq\gamma(\{z\})\subseteq\cl x$.
\end{definition}

\begin{lemma}[{\cite[Propositions 3.9--3.12 and 3.14, Theorems 3.13, 3.17, 4.7 and 4.8]{wangli}}]\label{lem:wang}
Let $X$ be an IG-closure space.
\begin{enumerate}
\item For every $x\in X$, the set $\cl x$ belongs to $\G(X)$. For $U\in\G(X)$, the family $\set{\cl x}{x\in U}$ is directed and $U={\bigcup}_{x\in U}\cl x$.
\item The poset $\G(X)$ is a domain; directed suprema are unions, and $\mathcal B_X=\set{\cl x}{x\in X}$ is a basis. For $U,V\in\G(X)$,
\[
U\ll V\quad\Longleftrightarrow\quad U\subseteq\cl x\text{ for some }x\in V.
\]
\item If $X$ is locally consistent, then $\G(X)$ is an $L$-domain.
\item Let $D$ be a domain. For $A\subseteq D$, put
\[
\cl A=\set{d\in D}{d\ll a\text{ for some }a\in A}.
\]
Then $(D,\cl\cdot)$ is an IG-closure space, and the map
\begin{equation}\label{eq:iso}
\eta_D\colon D\longrightarrow\G(D),\qquad d\longmapsto\dda d
\end{equation}
is an order isomorphism. If $D$ is an $L$-domain, this IG-closure space is locally consistent.
\end{enumerate}
\end{lemma}
Let $X$ be an IG-closure space and $D$ a domain. We call
$\mathcal B_X=\set{\cl x}{x\in X}$ the \emph{canonical basis} of
$\G(X)$, and we call the IG-closure space $(D,\cl\cdot)$
constructed in Lemma~\ref{lem:wang}(4) the \emph{canonical
IG-closure space} of $D$. For this space we have $\cl d=\dda d$ for every $d\in D$. Since $\eta_D$ in \eqref{eq:iso} is an order isomorphism onto $\G(D)$, we have $d\le e$ if and only if $\cl d\subseteq\cl e$ for all $d,e\in D$, and
\[
\mathcal B_D=\set{\cl d}{d\in D}=\G(D).
\]
\section{A characterization of Lawson compactness}\label{sec:lc}
Throughout this section, we assume that $X$ is an IG-closure space.
To formulate the conditions below, we define a relation $\triangleleft$
on $X$ by
\begin{equation}\label{eq:triangle}
x\triangleleft y\quad\Longleftrightarrow\quad \cl x\subseteq\cl t\text{ for some }t\in\cl y.
\end{equation}
By Lemma~\ref{lem:wang}(2), $x\triangleleft y$ is equivalent to $\cl x\ll\cl y$ in $\G(X)$. In particular, $x\triangleleft y$ implies $\cl x\subseteq\cl y$. 
\begin{lemma}\label{lem:finitegen}
Let $X$ be an IG-closure space. Then $\G(X)$ is finitely generated if and only if $X$ satisfies
\condition{(C1)}{there exists $M\fin X$ such that, for every $x\in X$, there exists $m\in M$ with $\cl m\subseteq\cl x$.}
\end{lemma}
\begin{proof}
Suppose that $\G(X)$ is finitely generated. Then there exist $U_1,\ldots,U_n\in\G(X)$ such that
\[
\G(X)=\ua_{\G(X)}\{U_1,\ldots,U_n\}.
\]
It follows from Lemma~\ref{lem:gd}(1) that each $U_i$ is nonempty, and we may choose $m_i\in U_i$.  By Lemma~\ref{lem:gd}(2), we have $\cl{m_i}\subseteq U_i$
for each $i\in\{1,\ldots,n\}$.
Now let $x\in X$. By Lemma~\ref{lem:wang}(1),
$\cl x\in\G(X)$, so $U_i\subseteq\cl x$ for some
$i\in\{1,\ldots,n\}$. This implies that
\[
\cl{m_i}\subseteq U_i\subseteq\cl x.
\]
Therefore, $X$ satisfies \textup{(C1)} with
$M=\{m_1,\ldots,m_n\}$.

Conversely, suppose that $X$ satisfies \textup{(C1)}.
Let $M\fin X$ be as in \textup{(C1)}.
By Lemma~\ref{lem:wang}(1), $\cl m\in\G(X)$ for every
$m\in M$.
Given $U\in\G(X)$, by Lemma~\ref{lem:gd}(1), we may choose
$x\in U$. Lemma~\ref{lem:gd}(2) gives $\cl x\subseteq U$.
By \textup{(C1)}, there exists $m\in M$ such that
\[
\cl m\subseteq\cl x\subseteq U.
\]
Consequently,
\[
\G(X)=\ua_{\G(X)}\set{\cl m}{m\in M}.
\]
Thus $\G(X)$ is finitely generated.
\end{proof}
\begin{lemma}\label{lem:M}
Let $X$ be an IG-closure space. Then $\G(X)$ satisfies the second condition of Lemma~\ref{lem:lawson} with respect to the canonical basis $\mathcal B_X$ if and only if $X$ satisfies
\condition{(C2)}{whenever $x_i,y_i\in X$ with $y_i\triangleleft x_i$ $(i=1,2)$, there exists $F\fin X$ such that}
\subcondition{(i)}{$\cl{y_1}\cup\cl{y_2}\subseteq\cl w$ for every $w\in F$;}
\subcondition{(ii)}{for every $z\in X$ with $x_i\triangleleft z$ for $i=1,2$, there exists $w\in F$ such that $\cl w\subseteq\cl z$.}
\end{lemma}
\begin{proof}
Throughout this proof, $\ua$ is taken in $\G(X)$.

\medskip\indent($\boldsymbol{\Rightarrow}$) Suppose that $\G(X)$ satisfies the second condition of Lemma~\ref{lem:lawson} with respect to $\mathcal B_X$. We show that $X$ satisfies \textup{(C2)}. Let $x_i,y_i\in X$ satisfy $y_i\triangleleft x_i$ for $i=1,2$. By the equivalence following \eqref{eq:triangle}, we have
\[
\cl{y_i}\ll\cl{x_i}\qquad(i=1,2).
\]
Since $\mathcal B_X=\set{\cl x}{x\in X}$, that condition yields $F\fin X$ such that
\begin{equation}\label{eq:M-witness}
\ua\cl{x_1}\cap\ua\cl{x_2}
\subseteq
\ua\set{\cl w}{w\in F}
\subseteq
\ua\cl{y_1}\cap\ua\cl{y_2}.
\end{equation}
For each $w\in F$, the second inclusion in \eqref{eq:M-witness} implies that $\cl w\in\ua\cl{y_1}\cap\ua\cl{y_2}$. Hence
\[
\cl{y_1}\cup\cl{y_2}\subseteq\cl w,
\]
so condition \textup{(C2)(i)} holds.
Now let $z\in X$ satisfy $\cl{x_i}\subseteq\cl z$ for $i=1,2$. Then
\[
\cl z\in\ua\cl{x_1}\cap\ua\cl{x_2}.
\]
By the first inclusion in \eqref{eq:M-witness}, there exists $w\in F$ such that $\cl w\subseteq\cl z$. In particular, this conclusion holds whenever $x_i\triangleleft z$ for $i=1,2$, because $x_i\triangleleft z$ implies $\cl{x_i}\subseteq\cl z$. Thus condition \textup{(C2)(ii)} also holds, and $X$ satisfies \textup{(C2)}.

\medskip\indent($\boldsymbol{\Leftarrow}$) Suppose that $X$ satisfies \textup{(C2)}.   Let $x_i,y_i\in X$ satisfy
\[
\cl{y_i}\ll\cl{x_i}\qquad(i=1,2).
\]
Since $\mathcal B_X$ is a basis of $\G(X)$, basis interpolation yields $p_i\in X$ such that
\[
\cl{y_i}\ll\cl{p_i}\ll\cl{x_i}
\qquad(i=1,2).
\]
By the equivalence following \eqref{eq:triangle}, we have $y_i\triangleleft p_i$ for $i=1,2$. Applying \textup{(C2)} with $p_i$ in place of $x_i$ $(i=1,2)$, we obtain $F\fin X$ satisfying \textup{(C2)(i)} and \textup{(C2)(ii)}; in particular, $\cl{y_1}\cup\cl{y_2}\subseteq\cl w$ for every $w\in F$.
By \textup{(C2)(i)}, each $w\in F$ satisfies $\cl w\in\ua\cl{y_1}\cap\ua\cl{y_2}$. Since the latter set is an upper set,
\begin{equation}\label{eq:lowerincl}
	\ua\set{\cl w}{w\in F}
	\subseteq
	\ua\cl{y_1}\cap\ua\cl{y_2}.
\end{equation}
It remains to prove the first inclusion in \eqref{eq:M-witness}
\begin{equation}\label{eq:rev}
	\ua\cl{x_1}\cap\ua\cl{x_2}\subseteq\ua\set{\cl w}{w\in F}.
\end{equation}
To this end, let $U\in\ua\cl{x_1}\cap\ua\cl{x_2}$. Then
$\cl{x_i}\subseteq U$ for $i=1,2$. Since
$\cl{p_i}\ll\cl{x_i}$, it follows that
$\cl{p_i}\ll U$. By Lemma~\ref{lem:wang}(2), there
exists $z_i\in U$ such that
\[
\cl{p_i}\subseteq\cl{z_i}
\qquad(i=1,2).
\]
As $\{z_1,z_2\}\fin U$ and $U$ is generalized directed,
there exists $z\in U$ such that
\[
\{z_1,z_2\}\subseteq\cl z\subseteq U.
\]
For each $i=1,2$, we have
$z_i\in\cl z$ and
$\cl{p_i}\subseteq\cl{z_i}$. Thus
$p_i\triangleleft z$ by \eqref{eq:triangle}. Condition \textup{(C2)(ii)} therefore gives an element $w\in F$ such that
\[
\cl w\subseteq\cl z\subseteq U.
\]
Consequently, $U\in\ua\set{\cl w}{w\in F}$. Since $U$ was arbitrary, \eqref{eq:rev} holds.
Together with \eqref{eq:lowerincl}, this yields \eqref{eq:M-witness}, and hence $\G(X)$ satisfies the second condition of Lemma~\ref{lem:lawson} with respect to $\mathcal B_X$.
\end{proof}
The following examples show that conditions \textup{(C1)} and \textup{(C2)} are independent.

\begin{example}\label{ex:C1notC2}
	Let
	\[
	X_1=\{p,a,b\}\cup\set{c_n}{n\in\mathbb N},
	\]
	where $p$, $a$, $b$, and the $c_n$ are pairwise distinct.
	Define an operation $\cl\cdot:\mathcal P(X_1)\to\mathcal P(X_1)$
	by
	\[
	\cl A={\bigcup}_{x\in A}C_x
	\qquad(A\subseteq X_1),
	\]
	where
	\[
	C_p=\{p\},\qquad
	C_a=\{p,a\},\qquad
	C_b=\{p,b\},\qquad
	C_{c_n}=\{p,a,b,c_n\}
	\quad(n\in\mathbb N).
	\]
	Then $\cl x=C_x$ for every $x\in X_1$. The operation
	$\cl\cdot$ is monotone and idempotent, and $x\in\cl x$ for every
	$x\in X_1$. If $x\in X_1$ and $M\fin\cl x$, then the choice $y=x$
	gives $y\in\cl x$ and $M\subseteq\cl y$. Thus \textup{(In)}
	holds, and $X_1$ is an IG-closure space.
	
	Since $\cl p=\{p\}\subseteq\cl x$ for every $x\in X_1$,
	condition \textup{(C1)} holds with $M=\{p\}$.
	
	We show that \textup{(C2)} fails. By \eqref{eq:triangle},
	\[
	a\triangleleft a,\qquad b\triangleleft b,\qquad
	a\triangleleft c_n,\qquad b\triangleleft c_n
	\quad(n\in\mathbb N).
	\]
	Take $x_1=y_1=a$ and $x_2=y_2=b$, and suppose that
	$F\fin X_1$ satisfies \textup{(C2)(i)} and \textup{(C2)(ii)}.
	By \textup{(C2)(i)}, every $w\in F$ satisfies
	\[
	\{p,a,b\}=\cl a\cup\cl b\subseteq\cl w.
	\]
	Hence every $w\in F$ is of the form $c_k$ for some
	$k\in\mathbb N$. For each $n\in\mathbb N$, condition
	\textup{(C2)(ii)} applied to $z=c_n$ gives a $w\in F$ such that
	$\cl w\subseteq\cl{c_n}$. Since $w=c_k$ for some
	$k\in\mathbb N$, this inclusion forces $k=n$. Thus $c_n\in F$
	for every $n\in\mathbb N$, contradicting the finiteness of $F$.
	Therefore $X_1$ does not satisfy \textup{(C2)}.
\end{example}

\begin{example}\label{ex:C2notC1}
	Let $X_2$ be an infinite set and define
	\[
	\cl A=A
	\qquad(A\subseteq X_2).
	\]
	The operation is monotone and idempotent, and $x\in\cl x$ for
	every $x\in X_2$. Thus \textup{(In)} holds with $y=x$, and
	$X_2$ is an IG-closure space. Moreover,
	\[
	\cl x=\{x\}
	\qquad\text{and}\qquad
	x\triangleleft y\iff x=y
	\quad(x,y\in X_2).
	\]
	
	To verify \textup{(C2)}, let $x_i,y_i\in X_2$ satisfy
	$y_i\triangleleft x_i$ for $i=1,2$. Then $y_i=x_i$ for
	$i=1,2$. If $x_1=x_2$, the set $F=\{x_1\}$ satisfies
	\textup{(C2)(i)} and \textup{(C2)(ii)}. If $x_1\ne x_2$,
	there is no $z\in X_2$ such that $x_i\triangleleft z$ for
	$i=1,2$, so $F=\varnothing$ satisfies both
	\textup{(C2)(i)} and \textup{(C2)(ii)}.
	
	On the other hand, \textup{(C1)} fails. Since $\cl x=\{x\}$
	for every $x\in X_2$, any set $M\subseteq X_2$ satisfying the
	requirement in \textup{(C1)} must contain every element of
	$X_2$. As $X_2$ is infinite, no such finite set exists.
\end{example}
\begin{remark}\label{rem:lc2}
Let \textup{(C2$'$)} be the condition obtained from \textup{(C2)} by replacing the requirements $x_1\triangleleft z$ and $x_2\triangleleft z$ in \textup{(C2)(ii)} with
\[
\cl{x_1}\cup\cl{x_2}\subseteq\cl z.
\]
Since $x_i\triangleleft z$ implies $\cl{x_i}\subseteq\cl z$ for $i=1,2$, condition \textup{(C2$'$)} implies \textup{(C2)}.

Conversely, suppose that $X$ satisfies \textup{(C2)}. By Lemma~\ref{lem:M}, $\G(X)$ satisfies the second condition of Lemma~\ref{lem:lawson} with respect to $\mathcal B_X$. The derivation of \textup{(C2)} from that condition in the proof of Lemma~\ref{lem:M} also establishes \textup{(C2$'$)}. Thus \textup{(C2)} and \textup{(C2$'$)} are equivalent.

\end{remark}
\begin{theorem}\label{thm:main}
Let $X$ be an IG-closure space. Then $\G(X)$ is Lawson compact if and only if $X$ satisfies \textup{(C1)} and \textup{(C2)}.
\end{theorem}

\begin{proof}
By Lemma~\ref{lem:wang}(2), $\G(X)$ is a domain with basis
$\mathcal B_X$. The conclusion therefore follows from
Lemmas~\ref{lem:lawson}, \ref{lem:finitegen}
and~\ref{lem:M}.
\end{proof}

\begin{theorem}\label{thm:repr_domain}
A poset $D$ is a Lawson compact domain if and only if it is order-isomorphic to $\G(X)$ for an IG-closure space $X$ satisfying \textup{(C1)} and \textup{(C2)}.
\end{theorem}
\begin{proof}
Suppose that $D$ is a Lawson compact domain, and let $X$ be its canonical IG-closure space. By Lemma~\ref{lem:wang}(4), $\eta_D\colon D\to\G(X)$ is an order isomorphism. It is a homeomorphism for the Lawson topologies: it preserves Scott-open sets and principal filters, as does its inverse. Hence $\G(X)$ is Lawson compact, and Theorem~\ref{thm:main} gives \textup{(C1)} and \textup{(C2)}. 

Conversely, if $D$ is order-isomorphic to $\G(X)$ for an IG-closure space $X$ satisfying these conditions, then $\G(X)$ is a Lawson compact domain by Lemma~\ref{lem:wang}(2) and Theorem~\ref{thm:main}. Thus $D$ is a Lawson compact domain.
\end{proof}

\begin{example}\label{ex:unit}
Let $D=[0,1]$ with its usual order. Then $D$ is a continuous lattice and hence an $L$-domain. Its Lawson topology is the usual topology on $[0,1]$, so $D$ is Lawson compact.

Let $X$ be the canonical IG-closure space of $D$. Then
\[
\cl 0=\{0\},
\qquad
\cl x=[0,x)\quad(x>0),
\]
and
\[
x\triangleleft y
\quad\Longleftrightarrow\quad
x<y\ \text{or}\ x=0.
\]
Condition \textup{(C1)} holds with $M=\{0\}$.

To verify \textup{(C2)}, let $x_i,y_i\in D$ satisfy $y_i\triangleleft x_i$ for $i=1,2$. Put
\[
w=\max\{y_1,y_2\},
\qquad
F=\{w\}.
\]
Since $y_1,y_2\le w$, we have $\cl{y_1}\cup\cl{y_2}\subseteq\cl w$, so \textup{(C2)(i)} holds. Now let $z\in D$ satisfy $x_i\triangleleft z$ for $i=1,2$. Then $y_i\le x_i\le z$ for both $i$, and hence $w\le z$. Therefore $\cl w\subseteq\cl z$, proving \textup{(C2)(ii)}.
\end{example}

\section{Bases and minimal upper bounds}\label{sec:m}
In this section we answer \textup{(Q3)}, and \textup{(Q4)} as a consequence (Theorem~\ref{thm:selected}). We begin by recalling properties $m$ and $M$ for posets.

Let $P$ be a poset. For $A\fin P$, write
\[
\UB_{P}(A)=\set{p\in P}{a\leq p\text{ for all }a\in A},
\]
and write $\mub_{P}(A)$ for the set of minimal elements of $\UB_{P}(A)$.
\begin{definition}[{\cite[Section~1.4]{jung}}]\label{def:property-m}
A poset $P$ has \emph{property $m$} if, for every $A\fin P$ and every $z\in\UB_P(A)$, there exists $m\in\mub_P(A)$ such that $m\leq z$.
The poset $P$ has \emph{property $M$} if it has property $m$ and
$\mub_P(A)$ is finite for every $A\fin P$; see also
\cite[Definition~2.6]{jiang}.
\end{definition}

The following elementary observation will be useful.
\begin{lemma}\label{lem:mubfinite}
A poset $P$ has property $M$ if and only if, for every $A\fin P$, there exists $H\fin P$ such that
\begin{equation}\label{eq:ubgen}
\UB_{P}(A)=\ua_P H.
\end{equation}
\end{lemma}
\begin{proof}
Suppose that $P$ has property $m$ and that $\mub_P(A)$ is finite for every $A\fin P$. For such a set $A$, property $m$ gives
$\UB_P(A)=\ua_P\mub_P(A)$.
Thus \eqref{eq:ubgen} holds with $H=\mub_P(A)$.

Conversely, let $A\fin P$ and suppose that \eqref{eq:ubgen} holds for some $H\fin P$. Denote by $\operatorname{Min}(H)$  the set of minimal elements of $H$. Since $H$ is finite,
\[
\ua_P H=\ua_P\operatorname{Min}(H).
\]
Moreover, the minimal elements of $\ua_P H$ are precisely the elements of $\operatorname{Min}(H)$. Consequently,
\[
\mub_P(A)=\operatorname{Min}(H),
\qquad
\UB_P(A)=\ua_P\mub_P(A).
\]
The first equality shows that $\mub_P(A)$ is finite, while the second shows that every upper bound of $A$ lies above a minimal upper bound of $A$. Since $A$ was arbitrary, $P$ has property $m$ and $\mub_P(A)$ is finite for every $A\fin P$.
\end{proof}
The next lemma is the link with Lawson compactness.
\begin{lemma}[{\cite[Remark~2.8]{jiang}}]\label{lem:fmlawson}
Let $D$ be a domain. If $D$ has a basis $\mathcal B$ with property $M$, then $D$ is Lawson compact.
\end{lemma}

\subsection{The canonical basis}
We first characterize when the canonical basis $\mathcal B_X$ of $\G(X)$ has property $M$. To this end, we introduce the following condition on an IG-closure space $X$.
\condition{(CM)}{For every $M\fin X$, there exists $H\fin X$ such that:}
\subcondition{(i)}{$\cl m\subseteq\cl w$ for every $m\in M$ and every $w\in H$;}
\subcondition{(ii)}{for every $z\in X$ with $\cl m\subseteq\cl z$ for all $m\in M$, there exists $w\in H$ such that $\cl w\subseteq\cl z$.}

\begin{proposition}\label{prop:cm}
Let $X$ be an IG-closure space. Then $X$ satisfies \textup{(CM)} if and only if the canonical basis $\mathcal B_X=\set{\cl x}{x\in X}$, ordered by inclusion, has property $M$.
\end{proposition}
\begin{proof}
For $M,H\fin X$, put
\[
A_M=\set{\cl m}{m\in M}
\qquad\text{and}\qquad
A_H=\set{\cl w}{w\in H}.
\]
Then \textup{(CM)(i)} holds if and only if
\[
\ua_{\mathcal B_X}A_H\subseteq\UB_{\mathcal B_X}(A_M),
\]
and \textup{(CM)(ii)} holds if and only if
\[
\UB_{\mathcal B_X}(A_M)\subseteq\ua_{\mathcal B_X}A_H.
\]
Thus $H$ satisfies \textup{(CM)(i)} and \textup{(CM)(ii)} for $M$ if and only if
\[
\UB_{\mathcal B_X}(A_M)=\ua_{\mathcal B_X}A_H.
\]

Every finite subset of $\mathcal B_X$ has the form $A_M$ for
some $M\fin X$. The stated equivalence therefore follows from
Lemma~\ref{lem:mubfinite}, applied to $\mathcal B_X$ ordered
by inclusion.
\end{proof}
\begin{corollary}\label{cor:CM}
Let $X$ be an IG-closure space satisfying \textup{(CM)}. Then $\G(X)$ is Lawson compact, and $X$ satisfies \textup{(C1)} and \textup{(C2)}.

\end{corollary}
\begin{proof}
By Lemma~\ref{lem:wang}(2), $\G(X)$ is a domain with canonical basis $\mathcal B_X$. Since $X$ satisfies \textup{(CM)}, Proposition~\ref{prop:cm} shows that $\mathcal B_X$ has property $M$. Hence $\G(X)$ is Lawson compact by Lemma~\ref{lem:fmlawson}. It follows from Theorem~\ref{thm:main} that $X$ satisfies \textup{(C1)} and \textup{(C2)}.
\end{proof}

\begin{theorem}\label{thm:whole}
	Let $D$ be a domain, and let $X$ be its canonical IG-closure
	space. Then $D$ has property $M$ if and only if $X$ satisfies
	\textup{(CM)}.
\end{theorem}

	\begin{proof}
		Since $X$ is the canonical IG-closure space of $D$,
		Lemma~\ref{lem:wang}(4) shows that the map
		$d\mapsto\cl d$ is an order isomorphism from $D$ onto
		$\mathcal B_X$. Therefore, $D$ has property $M$ if and only if
		$\mathcal B_X$ has property $M$. The assertion
		follows from Proposition~\ref{prop:cm}.
	\end{proof}

The next example shows that there exists an IG-closure
space $X$ satisfying \textup{(C1)} and \textup{(C2)} but
not \textup{(CM)}. Moreover, $X$ can be taken to be the
canonical IG-closure space of a Lawson compact algebraic
domain.

\begin{example}\label{ex:CMfails}
Let $N=\{0,1,2,\ldots\}\cup\{\omega\}$, where
$\{0,1,2,\ldots\}$ has its usual order and $\omega$ is the
greatest element. Write $\leq_N$ for the order on $N$, and let
$L=N\times\{0,1\}$ have the product order. Then $L$ is a
complete lattice. Set
\[
D=L\cup\set{c_n}{n\in N\setminus\{\omega\}},
\]
where the $c_n$ are pairwise distinct elements outside $L$.
The order $\leq_D$ on $D$ is illustrated in
Figure~\ref{fig:CMfails}: its restriction to $L$ is the
product order, each $c_n$ is maximal, and for every $k\in N$
and $n\in N\setminus\{\omega\}$,
\[
(k,0)\leq_D c_n,
\qquad
(k,1)\leq_D c_n\iff k\leq_N n.
\]

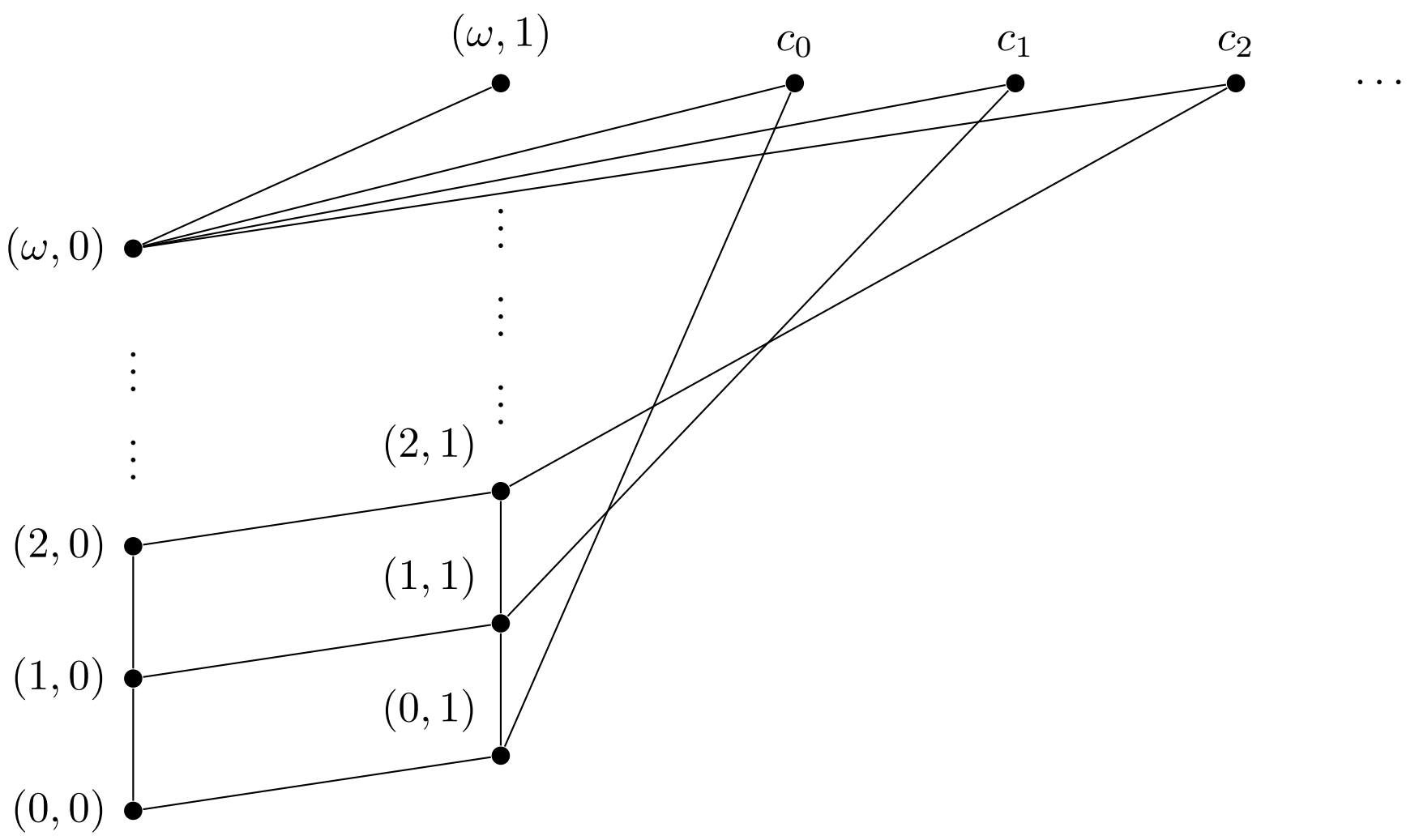
\begin{figure}[htbp]
	\centering
	\begin{tikzpicture}[
		x=1.2cm,
		y=0.9cm,
		every node/.style={font=\small},
		point/.style={circle, fill, inner sep=1.5pt}
		]
		
		\node[point, label=left:{$(0,0)$}] (a0) at (0,0) {};
		\node[point, label=left:{$(1,0)$}] (a1) at (0,1.2) {};
		\node[point, label=left:{$(2,0)$}] (a2) at (0,2.4) {};
		\node at (0,3.3) {$\vdots$};
		\node at (0,4.1) {$\vdots$};
		\node[point, label=left:{$(\omega,0)$}] (aw) at (0,5.1) {};
		
		\node[point, label distance=1pt, label=above left:{$(0,1)$}] (b0) at (2.5,0.5) {};
		\node[point, label distance=1pt, label=above left:{$(1,1)$}] (b1) at (2.5,1.7) {};
		\node[point, label distance=1pt, label=above left:{$(2,1)$}] (b2) at (2.5,2.9) {};
		\node at (2.5,3.8) {$\vdots$};
		\node at (2.5,4.6) {$\vdots$};
		\node at (2.5,5.4) {$\vdots$};
		
		\node[point, label=above:{$(\omega,1)$}] (bw) at (2.5,6.6) {};
		\node[point, label=above:{$c_0$}] (c0) at (4.5,6.6) {};
		\node[point, label=above:{$c_1$}] (c1) at (6,6.6) {};
		\node[point, label=above:{$c_2$}] (c2) at (7.5,6.6) {};
		\node at (8.5,6.6) {$\cdots$};
		
		\draw (a0) -- (a1) -- (a2);
		\draw (b0) -- (b1) -- (b2);
		
		\draw (a0) -- (b0);
		\draw (a1) -- (b1);
		\draw (a2) -- (b2);
		\draw (aw) -- (bw);
		
		\draw (aw) -- (c0);
		\draw (aw) -- (c1);
		\draw (aw) -- (c2);
		\draw (b0) -- (c0);
		\draw (b1) -- (c1);
		\draw (b2) -- (c2);
		
	\end{tikzpicture}
\caption{A schematic Hasse diagram of $D$.}
	\label{fig:CMfails}
\end{figure}

\noindent\textbf{Claim 1.} $D$ is an algebraic domain.

For each $n\in N\setminus\{\omega\}$, put
\[
L_n=\da_D c_n\cap L
=\da_L(\omega,0)\cup\da_L(n,1).
\]
Then each $L_n$ is a lower set of $L$ closed under directed suprema in $L$. A
directed subset of $D$ containing some $c_n$ has $c_n$ as its greatest element;
every other directed subset $E$ is contained in $L$. In the latter case, let
$v={\bigvee}_L E$. If $c_n$ is an upper bound of $E$, then $E\subseteq L_n$.
Since $L_n$ is closed under directed suprema in $L$, we have $v\in L_n$ and
hence $v\leq_D c_n$. Thus $v$ is also the supremum of $E$ in $D$. Therefore
$D$ is a dcpo.

The compact elements of $D$ are
\[
\K(D)=D\setminus\{(\omega,0),(\omega,1)\}.
\]
Moreover, for each $i\in\{0,1\}$, $(\omega,i)={\bigvee}_{k\in N\setminus\{\omega\}}(k,i)$.
Hence $D$ is algebraic.

\medskip
\noindent\textbf{Claim 2.} $D$ is Lawson compact.

We first show that, for all
$a_1,a_2\in\K(D)$, there exists $F\fin \K(D)$ such that
\[
\ua_D a_1\cap\ua_D a_2=\ua_D F.
\]
Up to interchanging $a_1$ and $a_2$, we consider the following
cases.
\begin{itemize}
\item For comparable $a_1$ and $a_2$, take $F$ to be the
singleton containing the larger element.
\item For $(k,0)$ and $(j,1)$, where
$j,k\in N\setminus\{\omega\}$ and $j<_N k$, take
\(
F=\{(k,1),c_j,\ldots,c_{k-1}\}.
\)
\item For $(j,1)$ and $c_n$, where
$j,n\in N\setminus\{\omega\}$ and $n<_N j$, or for $c_m$
and $c_n$, where $m,n\in N\setminus\{\omega\}$ and
$m\ne n$, there is no common upper bound; take
$F=\varnothing$.
\end{itemize}
Now let $a_i,b_i\in\K(D)$ satisfy $b_i\ll_D a_i$ for $i=1,2$. Let
$F\fin\K(D)$ be as above for $a_1$ and $a_2$. Since $b_i\leq_D a_i$ for $i=1,2$, we have
\[
\ua_D a_1\cap\ua_D a_2=\ua_D F\subseteq\ua_D b_1\cap\ua_D b_2.
\]
Thus $D$ satisfies the second condition of Lemma~\ref{lem:lawson} with respect to $\K(D)$, and since
$D=\ua_D(0,0)$, Lemma~\ref{lem:lawson} shows that $D$ is Lawson compact.

\medskip
\noindent\textbf{Claim 3.} $\K(D)$ has property $M$.

Put $K=\K(D)$. By Lemma~\ref{lem:mubfinite}, it suffices
to show that, for every $A\fin K$, there exists $H\fin K$
such that
\[
\UB_K(A)=\ua_K H.
\]
We prove this by induction on the number of elements
of $A$.

For $A=\varnothing$, take $H=\{(0,0)\}$, since $(0,0)$ is
the least element of $K$. For a singleton $A=\{a\}$, take
$H=\{a\}$.

Suppose that $\UB_K(A)=\ua_K H$ for some $H\fin K$, and
let $b\in K\setminus A$. By Claim~2, for each $h\in H$ there
exists $F_h\fin K$ such that
$\UB_K(\{h,b\})=\ua_K F_h$. Hence
\[
\UB_K(A\cup\{b\})
=\UB_K(A)\cap\ua_K b
={\bigcup}_{h\in H}\UB_K(\{h,b\})
=\ua_K\Bigl({\bigcup}_{h\in H}F_h\Bigr).
\]
Since $H$ and each $F_h$ are finite, their union
${\bigcup}_{h\in H}F_h$ is finite. This completes the induction
and proves the claim.

\medskip
\noindent\textbf{Claim 4.} \textup{(C1)} and \textup{(C2)}
do not imply \textup{(CM)}.

Let $X$ be the canonical IG-closure space of $D$. By
Lemma~\ref{lem:wang}(4), the map
\[
\eta_D\colon D\longrightarrow\mathcal B_X,\qquad
d\longmapsto\cl d,
\]
is an order isomorphism, and $\G(X)=\mathcal B_X$. Since $D$
is Lawson compact, so is $\G(X)$. Hence $X$ satisfies
\textup{(C1)} and \textup{(C2)} by
Theorem~\ref{thm:main}.

Now let $A=\{(\omega,0),(0,1)\}\subseteq D$. Its upper bounds in $D$
are $(\omega,1)$ and the elements $c_n$ with
$n\in N\setminus\{\omega\}$. These elements are pairwise
incomparable, so
\[
\mub_D(A)
=\{(\omega,1)\}\cup
\set{c_n}{n\in N\setminus\{\omega\}},
\]
which is infinite. Since $\eta_D$ is an order isomorphism, the
finite set $\eta_D(A)$ has infinitely many minimal upper bounds
in $\mathcal B_X$. By Proposition~\ref{prop:cm}, $X$ does not
satisfy \textup{(CM)}.\end{example}

\subsection{Bases determined by subsets of \texorpdfstring{$X$}{X}}
We now consider domains that admit a basis with property $M$. As Example~\ref{ex:CMfails} shows, the canonical basis $\mathcal B_X$ of $\G(X)$ need not have property $M$, even when another basis does. To describe such bases, we consider subsets $S\subseteq X$.

Let $X$ be an IG-closure space. For $S\subseteq X$, put
\begin{equation}\label{eq:BS}
\mathcal B_S=\set{\cl s}{s\in S}.
\end{equation}
\begin{lemma}\label{lem:subbasis}
Let $X$ be an IG-closure space and $S\subseteq X$. Then $\mathcal B_S$ (see \eqref{eq:BS}) is a basis of $\G(X)$ if and only if the following condition holds:
\condition{(B)}{For all $x,y\in X$ with $y\triangleleft x$, there exists $s\in S$ such that $y\triangleleft s\triangleleft x$.}
\end{lemma}

\begin{proof}
Suppose first that $\mathcal B_S$ is a basis of $\G(X)$, and
let $x,y\in X$ satisfy $y\triangleleft x$. By the equivalence
following \eqref{eq:triangle}, $\cl y\ll\cl x$. Basis
interpolation for $\mathcal B_S$ gives $s\in S$ such that
$\cl y\ll\cl s\ll\cl x$. Hence
$y\triangleleft s\triangleleft x$, and \textup{(B)} holds.

Conversely, suppose that \textup{(B)} holds, and put
$D=\G(X)$. We show that $\mathcal B_S$ is a basis of $D$.
Let $U\in D$ and let $F\fin D$ satisfy $F\ll U$.

Since $\mathcal B_X$ is a basis of $D$ by
Lemma~\ref{lem:wang}(2), applying basis interpolation twice
gives $x,y\in X$ such that
$F\ll\cl y\ll\cl x\ll U$.
By the equivalence following \eqref{eq:triangle}, we have
$y\triangleleft x$. Thus \textup{(B)} gives $s\in S$ with
$y\triangleleft s\triangleleft x$. It follows that
$F\ll\cl s\ll U$, and $\cl s\in\mathcal B_S\cap\dda U$.

Taking $F=\varnothing$ shows that $\mathcal B_S\cap\dda U$ is
nonempty. Taking $F=\{V_1,V_2\}$ with
$V_1,V_2\in\mathcal B_S\cap\dda U$ shows that
$\mathcal B_S\cap\dda U$ is directed.

We now show that $U=\bigvee(\mathcal B_S\cap\dda U)$. For
every $V\in\mathcal B_S\cap\dda U$, we have $V\subseteq U$.
Let $W$ be an upper bound of $\mathcal B_S\cap\dda U$. For each
$V\in\mathcal B_X\cap\dda U$, taking $F=\{V\}$ gives $s\in S$
such that $\cl s\in\mathcal B_S\cap\dda U$ and
$V\subseteq\cl s$; hence $V\subseteq W$. Thus $W$ is an upper
bound of $\mathcal B_X\cap\dda U$. Since $\mathcal B_X$ is a
basis of $D$, the set $\mathcal B_X\cap\dda U$ has supremum
$U$, so $U\subseteq W$. Therefore
$U=\bigvee(\mathcal B_S\cap\dda U)$, and $\mathcal B_S$ is a
basis of $D$.
\end{proof}
\begin{lemma}\label{lem:ms}
Let $X$ be an IG-closure space and $S\subseteq X$. Then $(\mathcal B_S,\subseteq)$ has property $M$ if and only if the following condition holds:
\condition{(M$_S$)}{For every $M\fin S$, there exists $H\fin S$ such that:}
\subcondition{(i)}{$\cl m\subseteq\cl w$ for every $m\in M$ and every $w\in H$;}
\subcondition{(ii)}{for every $z\in S$ with $\cl m\subseteq\cl z$ for every $m\in M$, there exists $w\in H$ such that $\cl w\subseteq\cl z$.}
\end{lemma}
\begin{proof}
For $M,H\fin S$, put $A_M=\set{\cl m}{m\in M}$ and
$A_H=\set{\cl w}{w\in H}$. Since $\UB_{\mathcal B_S}(A_M)$
is an upper set of $\mathcal B_S$, condition
\textup{(M$_S$)(i)} is equivalent to
$\ua_{\mathcal B_S}A_H\subseteq\UB_{\mathcal B_S}(A_M)$,
and condition \textup{(M$_S$)(ii)} is equivalent to the
reverse inclusion. Thus $H$ satisfies \textup{(M$_S$)(i)}
and \textup{(M$_S$)(ii)} for $M$ if and only if
$\UB_{\mathcal B_S}(A_M)=\ua_{\mathcal B_S}A_H$. Every
finite subset of $\mathcal B_S$ has the form $A_M$ for some
$M\fin S$. The
conclusion follows from Lemma~\ref{lem:mubfinite}, applied
to $\mathcal B_S$.
\end{proof}
\begin{theorem}\label{thm:selected}
A poset $P$ is a domain admitting a basis $\mathcal B$ with
property $M$ if and only if $P$ is order isomorphic to
$\G(X)$ for an IG-closure space $X$ with a subset
$S\subseteq X$ satisfying \textup{(B)} and \textup{(M$_S$)}.
\end{theorem}
\begin{proof}
($\boldsymbol{\Rightarrow}$) Let $\mathcal B$ be a basis of $P$ with
property $M$. Let $X$ be the canonical IG-closure space of $P$,
and put $S=\mathcal B$. By Lemma~\ref{lem:wang}(4), the map $\eta_P$
in \eqref{eq:iso} is an order isomorphism from $P$ onto $\G(X)$, and
$\eta_P(\mathcal B)=\mathcal B_S$. Hence $\mathcal B_S$ is a basis of
$\G(X)$ with property $M$. By Lemmas~\ref{lem:subbasis}
and~\ref{lem:ms}, $S$ satisfies \textup{(B)} and \textup{(M$_S$)}.

\medskip\indent($\boldsymbol{\Leftarrow}$) Let $X$ be an IG-closure
space with a subset $S\subseteq X$ satisfying \textup{(B)} and
\textup{(M$_S$)}, and suppose that $P$ is order-isomorphic to $\G(X)$.
By Lemma~\ref{lem:wang}(2), $\G(X)$ is a domain. By
Lemmas~\ref{lem:subbasis} and~\ref{lem:ms}, $\mathcal B_S$ is a basis
of $\G(X)$ with property $M$. Since $P$ is order-isomorphic to
$\G(X)$, it is a domain admitting a basis with property $M$.
\end{proof}

\begin{remark}\label{rem:classes}
For an IG-closure space $X$, the conditions introduced above
are related by the following implications.
	\[
	\textup{(CM)}
	\ \Longrightarrow\
	\text{there exists }S\subseteq X\text{ satisfying
		\textup{(B)} and \textup{(M$_S$)}}
	\ \Longrightarrow\
	\textup{(C1) and (C2)}.
	\]
	
	The first implication follows by taking $S=X$. Indeed,
	\textup{(B)} holds by basis interpolation, and
	\textup{(M$_X$)} is precisely \textup{(CM)}. For the
	second implication, let $S\subseteq X$ satisfy
	\textup{(B)} and \textup{(M$_S$)}. By
	Lemmas~\ref{lem:subbasis} and~\ref{lem:ms},
	$\mathcal B_S$ is a basis of $\G(X)$ with property $M$. Thus $\G(X)$ is Lawson compact by
	Lemma~\ref{lem:fmlawson}, and $X$ satisfies
	\textup{(C1)} and \textup{(C2)} by
	Theorem~\ref{thm:main}.
	
	The first implication cannot be reversed in general.
	Let $D$ be the domain in Example~\ref{ex:CMfails}, and
	let $X$ be its canonical IG-closure space. That example
	shows that $X$ does not satisfy \textup{(CM)}, while
	$\K(D)$ is a basis of $D$ with property $M$.
	Put $S=\K(D)$. Since the order isomorphism
	$d\mapsto\cl d$ from $D$ onto $\G(X)$ maps $\K(D)$
	onto $\mathcal B_S$, the family $\mathcal B_S$ is a
	basis of $\G(X)$ and has property $M$. Hence $S$
	satisfies \textup{(B)} and \textup{(M$_S$)} by
	Lemmas~\ref{lem:subbasis} and~\ref{lem:ms}.
	
	The following example shows that the second implication
	cannot be reversed in general.
\end{remark}

\begin{example}\label{ex:nobasis}
Let $L=[0,1]\times\{0,1\}$ be equipped with the product
order $\leq_L$ induced by the usual order on $\mathbb R$.
Then $L$ is a complete lattice. Denote by $\mathcal I$ the set of
intervals $I=(p_I,q_I]$, where
\[
p_I=2^{-k}j,
\qquad
q_I=2^{-k}(j+1)
\]
for integers $k\geq 0$ and $0\leq j<2^k$. Set
\[
D=L\cup\set{c_I}{I\in\mathcal I},
\]
where the $c_I$ are pairwise distinct elements outside $L$.
Define an order $\leq_D$ on $D$ as follows:
\begin{itemize}
	\item for $(s,i),(t,j)\in L$,
	\[
	(s,i)\leq_D(t,j)
	\iff s\leq t\ \text{and}\ i\leq j;
	\]
	\item for $s\in[0,1]$ and $I\in\mathcal I$,
	\[
	(s,0)\leq_D c_I\iff s\leq q_I,
	\qquad
	(s,1)\leq_D c_I\iff s\leq p_I;
	\]
\item each $c_I$ is a maximal element of $D$.
\end{itemize}

To illustrate this order, fix $I\in\mathcal I$ with
$0<p_I<q_I<1$, and put
\[
P_I
=
\set{(s,i)}{s\in\{0,p_I,q_I,1\},\ i\in\{0,1\}}
\cup\{c_I\}.
\]
Figure~\ref{fig:nobasis} shows the order on $P_I$ inherited
from $D$.

\begin{figure}[htbp]
\centering
\begin{tikzpicture}[
x=1.4cm,
y=1.0cm,
every node/.style={font=\small},
point/.style={circle, fill, inner sep=1.5pt}
]

\node[point, label=left:{$(0,0)$}] (a0) at (0,0) {};
\node[point, label=left:{$(p_I,0)$}] (ap) at (0,1.4) {};
\node[point, label=left:{$(q_I,0)$}] (aq) at (0,2.8) {};
\node[point, label=left:{$(1,0)$}] (a1) at (0,4.2) {};

\node[point, label=right:{$(0,1)$}] (b0) at (3,0.9) {};
\node[point, label=right:{$(p_I,1)$}] (bp) at (3,2.3) {};
\node[point, label=right:{$(q_I,1)$}] (bq) at (3,3.7) {};
\node[point, label=above:{$(1,1)$}] (b1) at (3,5.1) {};

\node[point, label=above:{$c_I$}] (c) at (1.5,5.1) {};

\draw (a0) -- (ap) -- (aq) -- (a1);
\draw (b0) -- (bp) -- (bq) -- (b1);

\draw (a0) -- (b0);
\draw (ap) -- (bp);
\draw (aq) -- (bq);
\draw (a1) -- (b1);

\draw[preaction={draw=white, line width=3pt}]
(aq) -- (c);
\draw[preaction={draw=white, line width=3pt}]
(bp) -- (c);

\end{tikzpicture}
\caption{The order on $P_I$ inherited from $D$.}
\label{fig:nobasis}
\end{figure}
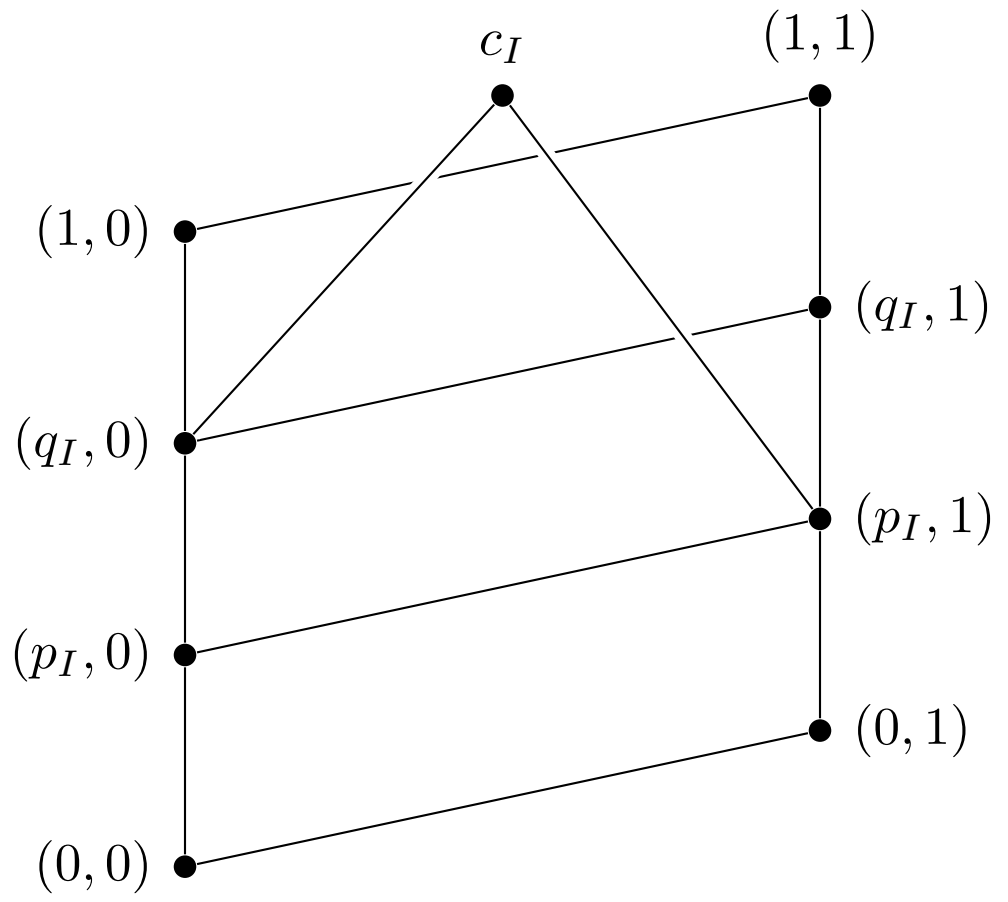

For each $I\in\mathcal I$, put
\[
L_I=\da_D c_I\cap L
=\da_L(q_I,0)\cup\da_L(p_I,1).
\]
Then $L_I$ is a lower set of $L$ and is closed under directed suprema in $L$. Indeed, the union of two lower sets is a lower set. Let $E\subseteq L_I$ be directed. Since $L$ is a complete
lattice, the supremum ${\bigvee}_L E$ exists. If $E\subseteq\da_L(q_I,0)$, then ${\bigvee}_L E\leq_L(q_I,0)$. Suppose that $E\nsubseteq\da_L(q_I,0)$, and choose $x\in E\setminus\da_L(q_I,0)$. For each $y\in E$, there exists a common upper bound $z\in E$ of $x$ and $y$. Since $x\leq_L z$, we have $z\notin\da_L(q_I,0)$. Thus $z\in\da_L(p_I,1)$, and
\(
y\leq_L z\leq_L(p_I,1).
\)
As $y$ was arbitrary, $E\subseteq\da_L(p_I,1)$, so ${\bigvee}_L E\leq_L(p_I,1)$. In either case, ${\bigvee}_L E\in L_I$.

\medskip
\noindent\textbf{Claim 1.} $D$ is a pointed $L$-domain.

A directed subset of $D$ containing some $c_I$ has $c_I$
as its greatest element, since $c_I$ is maximal. Every other
directed subset $E$ is contained in $L$. Let
$v={\bigvee}_L E$. If $c_I$ is an upper bound of $E$, then $E\subseteq L_I$.
Since $L_I$ is closed under directed suprema in $L$, we have
$v={\bigvee}_L E\in L_I$. Hence $v\leq_D c_I$. Thus $v$ is also the
supremum of $E$ in $D$. Consequently, $D$ is a dcpo.

 For every $c_I$, let $E\subseteq D$
be directed with $c_I\leq_D{\bigvee}_D E$. Since $c_I$ is maximal,
${\bigvee}_D E=c_I$. If $E\subseteq L$, then ${\bigvee}_D E\in L$,
a contradiction. Thus $E$ contains some $c_J$. By maximality,
$c_J$ is the greatest element of $E$, so
$c_J={\bigvee}_D E=c_I$. Hence $c_I$ is a compact element.

A directed subset $E$ of $D$ either contains some $c_I$, in
which case $c_I$ is its greatest element, or is contained in
$L$, in which case ${\bigvee}_D E={\bigvee}_L E$. For
$x,y\in L$, the definition of
the way-below relation gives $x\ll_D y$ if and only if
$x\ll_L y$. Thus,
for all $(s,i),(s',i')\in L$,
\[
(s,i)\ll_D(s',i')
\iff
i\leq i'
\quad\text{and}\quad
(s<s'\ \text{or}\ s=0).
\]
Moreover, since each $c_I$ is compact,
\[
\ell\ll_D c_I\iff \ell\leq_D c_I
\qquad(\ell\in D).
\]
 If $x\in L$, then $\dda_D x=\dda_L x$, which is directed
with supremum $x$ because $L$ is a continuous lattice. If $x=c_I$,
then $\dda_D c_I=\da_D c_I$ has greatest element $c_I$.
Thus $D$ is a domain with least element $(0,0)$.

For each $\ell\in L$, the principal ideal
$\da_D\ell=\da_L\ell$ is a complete lattice. Let
$I\in\mathcal I$. We show that
$\da_D c_I=L_I\cup\{c_I\}$ is also a complete lattice.
For a subset $S\subseteq L_I$, put $v={\bigvee}_L S$.
If $v\in L_I$, then $v$ is the supremum of $S$ in
$\da_D c_I$. If $v\notin L_I$, then $S$ has no upper
bound in $L_I$, since $L_I$ is a lower set of $L$.
In this case, $c_I$ is the supremum of $S$ in
$\da_D c_I$. Every subset of $\da_D c_I$ containing
$c_I$ has supremum $c_I$. Thus every subset of
$\da_D c_I$ has a supremum, and $\da_D c_I$ is a
complete lattice. Consequently, $D$ is an $L$-domain.

\medskip
\noindent\textbf{Claim 2.} $D$ is Lawson compact.

We apply Lemma~\ref{lem:lawson} with $D$ itself as a basis.
Clearly, $D=\ua_D(0,0)$. Let $a_i,b_i\in D$ satisfy
$b_i\ll_D a_i$ for $i=1,2$. We show that there exists $F\fin D$ such that
\[
\ua_D a_1\cap\ua_D a_2
\subseteq\ua_D F
\subseteq\ua_D b_1\cap\ua_D b_2.
\]

\noindent\textbf{Case 1.} Some $a_i$ is of the form $c_I$.

Since $c_I$ is maximal, $\ua_D c_I=\{c_I\}$. Thus
$\ua_D a_1\cap\ua_D a_2$ is either empty or equal to
$\{c_I\}$. Taking
$F=\ua_D a_1\cap\ua_D a_2$
gives the required finite set.

\noindent\textbf{Case 2.} $a_1,a_2\in L$.

Since $b_i\leq_D a_i$,
we also have $b_1,b_2\in L$. Write
\[
a_i=(\alpha_i,\varepsilon_i),
\qquad
b_i=(\beta_i,\delta_i)
,
\]
where $\alpha_i,\beta_i\in[0,1]$ and
$\varepsilon_i,\delta_i\in\{0,1\}$ for $i=1,2$.
Put $b=b_1\vee_L b_2$ and
\[
F=\{b\}\cup
\set{c_I}{a_1\leq_D c_I,\ a_2\leq_D c_I,\
b\nleq_D c_I}.
\]
Every element of $F$ lies above both $b_1$ and $b_2$.
Now let $u\in\ua_D a_1\cap\ua_D a_2$. If $u\in L$, then
\[
b=b_1\vee_L b_2\leq_L a_1\vee_L a_2\leq_L u.
\]
If $u=c_I$ for some $I\in\mathcal I$, then either
$b\leq_D c_I$ or $c_I\in F$. In both cases, $u\in\ua_D F$.
Therefore
\[
\ua_D a_1\cap\ua_D a_2
\subseteq\ua_D F
\subseteq\ua_D b_1\cap\ua_D b_2.
\]

It remains to show that $F$ is finite. If
$F\setminus\{b\}=\varnothing$, then $F=\{b\}$ is finite.
Suppose that $F\setminus\{b\}\neq\varnothing$. Write
$b=(\beta,\delta)$. If $\delta=0$, then every
$c_J\in F\setminus\{b\}$ satisfies
\[
\beta\leq\max\{\alpha_1,\alpha_2\}\leq q_J,
\]
and hence $b\leq_D c_J$, a contradiction. Thus $\delta=1$.
Since $\delta=\max\{\delta_1,\delta_2\}$, at least one of
$\delta_1,\delta_2$ equals $1$. By exchanging the indices
if necessary, we may assume that $\delta_2=1$.

Now let $c_I\in F\setminus\{b\}$ be arbitrary. Since
$b_2\leq_D a_2$, we have $\varepsilon_2=1$.
Moreover, $b\nleq_D c_I$ means $\beta>p_I$, and
$a_2\leq_D c_I$ gives
\[
\beta_2\leq\alpha_2\leq p_I<\beta.
\]
As $\beta=\max\{\beta_1,\beta_2\}$, it follows that
$\beta=\beta_1>p_I\geq0$. Since $b_1\ll_D a_1$,
we have $\beta_1<\alpha_1$. If $\varepsilon_1=1$,
then $a_1\leq_D c_I$ gives $\alpha_1\leq p_I$,
a contradiction. Hence $\varepsilon_1=0$, and
$a_1\leq_D c_I$ gives $\alpha_1\leq q_I$. Thus
\[
p_I<\beta_1<\alpha_1\leq q_I.
\]

Now $q_I-p_I=2^{-k}$ for some integer $k\geq0$, and
the preceding inequalities give
\[
2^{-k}=q_I-p_I>\alpha_1-\beta_1>0.
\]
The  inequality \(
2^{-k}>\alpha_1-\beta_1
\) holds for only finitely many
integers $k\geq0$. For each such $k$, there are $2^k$
intervals in $\mathcal I$ of length $2^{-k}$. Hence
$F\setminus\{b\}$ is finite, and so is $F$. Thus $D$
satisfies the second condition of Lemma~\ref{lem:lawson} with respect to $D$. Since $D=\ua_D(0,0)$,
Lemma~\ref{lem:lawson} shows that $D$ is Lawson compact.

\medskip
\noindent\textbf{Claim 3.} For every basis $\mathcal B$ of $D$,
$\mub_{\mathcal B}(A)$ is infinite for some $A\fin \mathcal B$.

Let $\mathcal B$ be a basis of $D$. Since every basis contains all
compact elements, we have $(0,1)\in \mathcal B$ and $c_I\in \mathcal B$ for
every $I\in\mathcal I$. Moreover,
\[
\dda_D(1,0)=\set{(s,0)}{0\leq s<1}.
\]
As $\mathcal B\cap\dda_D(1,0)$ has supremum $(1,0)$, there exists
$(r,0)\in \mathcal B$ with $0<r<1$. Put
$A=\{(r,0),(0,1)\}$.
Let $I\in\mathcal I$ be any interval containing $r$. Then $p_I<r\leq q_I$,
so $c_I$ is an upper bound of $A$. 

We show that $c_I$ is a minimal upper bound of $A$ in $D$.
Let $x\in D$ satisfy $x<_D c_I$. Then $x\in L_I$, so
$x=(s,0)$ or $x=(s,1)$ for some $s\in[0,1]$. If
$x=(s,0)$, then $(0,1)\nleq_D x$. If $x=(s,1)$, then
$s\leq p_I<r$, and hence $(r,0)\nleq_D x$. In either case,
$x$ is not an upper bound of $A$. Therefore
$c_I\in\mub_D(A)$. Since $A\subseteq \mathcal B$ and $c_I\in \mathcal B$,
it follows that $c_I\in\mub_{\mathcal B}(A)$.

For each integer $k\geq0$, there exists a unique integer $j_k$
with $0\leq j_k<2^k$ such that
\[
r\in I_k:=\bigl(2^{-k}j_k,\,2^{-k}(j_k+1)\bigr].
\]
The intervals $I_k$ belong to $\mathcal I$ and are pairwise
distinct, since they have different lengths. By the preceding
argument, $c_{I_k}\in\mub_{\mathcal B}(A)$ for every $k\geq0$.
The elements $c_{I_k}$ are pairwise distinct, so
$\mub_{\mathcal B}(A)$ is infinite.

\medskip
\noindent\textbf{Claim 4.} \textup{(C1)} and \textup{(C2)}
do not imply the existence of $S\subseteq X$
satisfying \textup{(B)} and \textup{(M$_S$)}.

Let $X$ be the canonical IG-closure space of $D$. Then
$\G(X)$ is order-isomorphic to $D$ and is therefore Lawson
compact. By Theorem~\ref{thm:main}, $X$ satisfies
\textup{(C1)} and \textup{(C2)}.

Suppose that a subset $S\subseteq X$ satisfies \textup{(B)}
and \textup{(M$_S$)}. By Lemmas~\ref{lem:subbasis}
and~\ref{lem:ms}, it determines a basis of $\G(X)$ with
property $M$. The order isomorphism
from $\G(X)$ onto $D$ then gives a basis of $D$ with
property $M$, contradicting Claim~3. Thus no such subset
$S$ exists. 
\end{example}
\begin{remark} Example~\ref{ex:nobasis} shows that a Lawson compact pointed $L$-domain need not have a basis with property $M$. Thus \textup{(Q1)} has a negative answer. \end{remark}

\section{The algebraic case}\label{sec:alg}
In this section, we study Lawson compactness and the basis
properties considered above in the algebraic setting. We give
representations of Lawson compact algebraic domains by dense
IG-closure spaces and obtain simpler conditions in the
locally consistent case. We also investigate when a dense
representing space can be chosen to satisfy \textup{(CM)}.
We begin by recalling the following notion of density for
IG-closure spaces.

\begin{definition}[{\cite[Definition~3.18]{wangli}}]\label{def:dense}
An IG-closure space $X$ is \emph{dense} if, whenever $b\in\cl a$, there exists $c\in\cl a$ such that $\{b,c\}\subseteq\cl c$.
\end{definition}
For an IG-closure space $X$, put
\begin{equation}\label{eq:Xc}
X_{\mathrm c}=\set{x\in X}{x\in\cl x}.
\end{equation}

The following result follows from the description of the canonical basis in \cite[Proposition~3.12]{wangli} and the representation of algebraic domains in \cite[Theorem~3.20]{wangli}.

\begin{lemma}\label{lem:compacts}
For every IG-closure space $X$,
\[
\K(\G(X))=\set{\cl x}{x\in X_{\mathrm c}}.
\]
Moreover, $\G(X)$ is algebraic if and only if $X$ is dense.
\end{lemma}

For algebraic domains, Lawson compactness can be characterized in
terms of the poset of compact elements.
\begin{lemma}[{\cite[Corollary~4.19]{jung}}]\label{lem:jung-algebraic}
Let $D$ be an algebraic domain. Then $D$ is Lawson compact if and
only if $\K(D)$ has property $M$.
\end{lemma}

\begin{lemma}[{\cite[Corollary~4.19]{jung}}]\label{lem:algcrit}
	An algebraic domain $D$ is Lawson compact if and only if
	$D=\ua H$ for some $H\fin\K(D)$ and, for any
	$a,b\in\K(D)$, there exists $F\fin\K(D)$ such that
	\[
	\ua a\cap\ua b=\ua F.
	\]
\end{lemma}

 To characterize Lawson
compactness of algebraic domains, we consider the following
conditions on an IG-closure space $X$.
\condition{(A1)}{There exists $M\fin X_{\mathrm c}$ such that, for every $x\in X_{\mathrm c}$, there exists $m\in M$ with $\cl m\subseteq\cl x$.}
\condition{(A2)}{For every $x,y\in X_{\mathrm c}$, there exists $F\fin X_{\mathrm c}$ such that}
\subcondition{(i)}{$\cl x\cup\cl y\subseteq\cl w$ for every $w\in F$;}
\subcondition{(ii)}{for every $z\in X_{\mathrm c}$ with $\cl x\cup\cl y\subseteq\cl z$, there exists $w\in F$ such that $\cl w\subseteq\cl z$.}
\begin{theorem}\label{thm:algebraic}
Let $X$ be a dense IG-closure space. Then $\G(X)$ is Lawson compact if and only if $X$ satisfies \textup{(A1)} and \textup{(A2)}.
\end{theorem}
\begin{proof}
	Put $D=\G(X)$. Since $X$ is dense, Lemma~\ref{lem:compacts} shows that
	$D$ is algebraic and that
	\[
	\K(D)=\set{\cl x}{x\in X_{\mathrm c}}.
	\]

	\medskip\indent($\boldsymbol{\Rightarrow}$) Suppose that $D$ is Lawson compact.
	By Lemma~\ref{lem:algcrit}, there exists $H\fin\K(D)$ such that
	$D=\ua H$.
	For each $h\in H$, choose $m_h\in X_{\mathrm c}$ with
	$h=\cl{m_h}$, and put
	\[
	M=\set{m_h}{h\in H}.
	\]
	
	Let $x\in X_{\mathrm c}$. Since $\cl x\in D=\ua H$,
	there exists $h\in H$ such that
	$\cl{m_h}=h\subseteq\cl x$.
	Thus $X$ satisfies \textup{(A1)}.

	To verify \textup{(A2)}, let $x,y\in X_{\mathrm c}$.
	By Lemma~\ref{lem:algcrit}, there exists $\mathcal F\fin \K(D)$ such that
	\[
	\ua\cl x\cap\ua\cl y=\ua\mathcal F.
	\]
	For each $V\in\mathcal F$, choose $w_V\in X_{\mathrm c}$
	with $V=\cl{w_V}$, and put
	$F=\set{w_V}{V\in\mathcal F}$.
	For every $w\in F$, the equality above gives
	$\cl x\cup\cl y\subseteq\cl w$.
	Now let $z\in X_{\mathrm c}$ satisfy
	$\cl x\cup\cl y\subseteq\cl z$. Then
	$\cl z\in\ua\mathcal F$, so there exists $V\in\mathcal F$
	such that
	$\cl{w_V}=V\subseteq\cl z$.
	Therefore $X$ satisfies \textup{(A2)}.

	\medskip\indent($\boldsymbol{\Leftarrow}$) Suppose that $X$ satisfies \textup{(A1)}
	and \textup{(A2)}. We show that $D$ is Lawson compact.
	By \textup{(A1)}, there exists $M\fin X_{\mathrm c}$ such that,
	for every $x\in X_{\mathrm c}$, there exists $m\in M$ with
	$\cl m\subseteq\cl x$. Put $H=\set{\cl m}{m\in M}$. Then $H$ is
	a finite subset of $\K(D)$. Let $U\in D$.
	Since $D$ is algebraic, there exists a compact element below $U$,
	say $\cl x$ with $x\in X_{\mathrm c}$. By \textup{(A1)}, there
	exists $m\in M$ such that
	\[
	\cl m\subseteq\cl x\subseteq U.
	\]
	Consequently, $D=\ua H$.

	Now let $a,b\in\K(D)$. Choose $x,y\in X_{\mathrm c}$
	such that $a=\cl x$ and $b=\cl y$. By \textup{(A2)},
	there exists $F\fin X_{\mathrm c}$
	satisfying \textup{(A2)(i)} and \textup{(A2)(ii)}.
	Put
	\[
	\mathcal F=\set{\cl w}{w\in F}.
	\]
	Then $\mathcal F$ is a finite subset of $\K(D)$.
	By \textup{(A2)(i)},
	\[
	\ua\mathcal F\subseteq\ua\cl x\cap\ua\cl y.
	\]

	For the reverse inclusion, let
	$U\in\ua\cl x\cap\ua\cl y$. Since $D$ is algebraic,
	the compact elements below $U$ form a directed set.
	It contains $\cl x$ and $\cl y$, so there exists
	$k\in\K(D)$ such that
	\[
	\cl x\cup\cl y\subseteq k\subseteq U.
	\]
	Choose $z\in X_{\mathrm c}$ with $k=\cl z$.
	By \textup{(A2)(ii)}, there exists $w\in F$ such that
	\[
	\cl w\subseteq\cl z\subseteq U.
	\]
	Hence $U\in\ua\mathcal F$. We have therefore shown that
	$\ua a\cap\ua b=\ua\mathcal F$. Together with $D=\ua H$,
	this proves that $D$ is Lawson compact by
	Lemma~\ref{lem:algcrit}.
\end{proof}

Let $X$ be a dense IG-closure space. By
Theorem~\ref{thm:algebraic}, $X$ satisfies \textup{(A1)}
and \textup{(A2)} if and only if $\G(X)$ is Lawson compact.
By Theorem~\ref{thm:main}, the latter is equivalent to $X$
satisfying \textup{(C1)} and \textup{(C2)}.

\begin{corollary}\label{cor:algebraic_meet}
For a dense IG-closure space $X$, the space $X$ satisfies \textup{(A1)} and \textup{(A2)} if and only if $\G(X)$ admits a basis with property $M$. 
\end{corollary}

\begin{proof}
	By Lemma~\ref{lem:compacts}, $\G(X)$ is an algebraic domain.
	Theorem~\ref{thm:algebraic} and Lemma~\ref{lem:jung-algebraic}
	show that, if $X$ satisfies \textup{(A1)} and \textup{(A2)},
	then the basis $\K(\G(X))$ has property $M$.
	Conversely, if $\G(X)$ has a basis with property $M$,
	then it is Lawson compact by Lemma~\ref{lem:fmlawson}.
	Hence $X$ satisfies \textup{(A1)} and \textup{(A2)} by
	Theorem~\ref{thm:algebraic}.
\end{proof}

\begin{theorem}\label{thm:alg_domain}
A poset $D$ is a Lawson compact algebraic domain if and only if it is order-isomorphic to $\G(X)$ for a dense IG-closure space $X$ satisfying \textup{(A1)} and \textup{(A2)}.
\end{theorem}
\begin{proof}
If $X$ is dense and satisfies the two conditions, Lemma~\ref{lem:compacts} and Theorem~\ref{thm:algebraic} give the conclusion. Conversely, let $D$ be a Lawson compact algebraic domain and $X$ its canonical IG-closure space. Lemma~\ref{lem:wang}(4) shows that $D$ is order-isomorphic to $\G(X)$; Lemma~\ref{lem:compacts} then implies that $X$ is dense. Theorem~\ref{thm:algebraic} yields \textup{(A1)} and \textup{(A2)}.
\end{proof}

We now give a criterion for Lawson compactness that uses the local structure of an algebraic $L$-domain. Let $X$ be an IG-closure space and put
\[
\mathcal K_X=\set{\cl x}{x\in X_{\mathrm c}}.
\]
The set $\mathcal K_X$ is ordered by inclusion. For $M\fin X_{\mathrm c}$, put
\[
\mathcal M_X(M)=\mub_{\mathcal K_X}\bigl(\set{\cl m}{m\in M}\bigr).
\]
In particular, $\mathcal M_X(\varnothing)$ is the set of minimal elements of $\mathcal K_X$.

\begin{theorem}\label{prop:algL-mub}
Let $X$ be a dense, locally consistent IG-closure space. Then $\G(X)$ is Lawson compact if and only if:
\condition{(L1)}{$\mathcal M_X(\varnothing)$ is finite;}
\condition{(L2)}{$\mathcal M_X(\{x,y\})$ is finite for all $x,y\in X_{\mathrm c}$.}
\end{theorem}
\begin{proof}
	Put $D=\G(X)$ and $K=\K(D)$. Since $X$ is dense and locally
	consistent, Lemma~\ref{lem:compacts} and
	Lemma~\ref{lem:wang}(3) show that $D$ is an algebraic
	$L$-domain and that
	\[
	K=\mathcal K_X=\set{\cl x}{x\in X_{\mathrm c}}.
	\]

	\medskip\indent($\boldsymbol{\Rightarrow}$) Suppose that $D$ is Lawson compact. By
	Lemma~\ref{lem:jung-algebraic}, $\mub_K(A)$ is finite for
	every $A\fin K$. Taking $A=\varnothing$ gives \textup{(L1)},
	and taking $A=\{\cl x,\cl y\}$ for $x,y\in X_{\mathrm c}$
	gives \textup{(L2)}.

	\medskip\indent($\boldsymbol{\Leftarrow}$) Suppose that \textup{(L1)} and \textup{(L2)} hold.
	We first show that $K$ has property $m$.
	
	 For any $A\fin K$,
	let $z\in K$ be an upper bound of $A$. Since $\da_D z$ is a
	complete lattice, $A$ has a supremum $u$ in $\da_D z$.
	If $w\leq u$ is an upper bound of $A$ in $D$, then
	$w\in\da_D z$, so $u\leq w$. Thus $u$ is a minimal upper
	bound of $A$ in $D$.  Since $A$ and $z$ were
	arbitrary, it suffices to check that $u$ is compact.

	Let $E\subseteq D$ be directed,
	put $s=\bigvee E$, and suppose that $u\leq s$. Let $v$ be the
	supremum of $A$ in $\da_D s$. Since $u$ is an upper bound of
	$A$ in $\da_D s$, we have $v\leq u$. As $v$ is also an upper
	bound of $A$, the minimality of $u$ gives $v=u$.
	For each $a\in A$, compactness of $a$ and $a\leq\bigvee E$
	give an element $e_a\in E$ with $a\leq e_a$. Since $A$ is
	finite and $E$ is directed, we may choose $e\in E$ above
	every $e_a$. If $A=\varnothing$, take any $e\in E$.
	Then $e$ is an upper bound of $A$ in $\da_D s$, so
	$u=v\leq e$. Hence $u$ is compact.

	We next prove by induction on the number of elements of $A$
	that $\UB_K(A)$ is finitely generated for every $A\fin K$.
	For $A=\varnothing$, property $m$ gives
	\[
	\UB_K(\varnothing)=\ua_K\mub_K(\varnothing),
	\]
	and $\mub_K(\varnothing)$ is finite by \textup{(L1)}.

	Suppose that $\UB_K(A)=\ua_K H$ for some $H\fin K$, and let
	$b\in K\setminus A$. For each $h\in H$, the set
	$\mub_K(\{h,b\})$ is finite by \textup{(L2)}, since $h$ and
	$b$ are closures of points in $X_{\mathrm c}$. By property
	$m$, we have $\UB_K(\{h,b\})=\ua_K\mub_K(\{h,b\})$. Hence
	\[
	\begin{aligned}
	\UB_K(A\cup\{b\})
	&=\UB_K(A)\cap\ua_K b
	{=\bigcup}_{h\in H}\UB_K(\{h,b\})\\
	&=\ua_K\Bigl({\bigcup}_{h\in H}\mub_K(\{h,b\})\Bigr).
	\end{aligned}
	\]
	The set ${\bigcup}_{h\in H}\mub_K(\{h,b\})$ is finite.
	This completes the induction.

	By Lemma~\ref{lem:mubfinite}, $\mub_K(A)$ is finite for
	every $A\fin K$. Since $K$ also has property $m$,
	Lemma~\ref{lem:jung-algebraic} shows that $D$ is Lawson
	compact.
\end{proof}

\begin{theorem}\label{thm:algL-mub-repr}
A poset $D$ is a Lawson compact algebraic $L$-domain if and only if it is order-isomorphic to $\G(X)$ for a dense, locally consistent IG-closure space $X$ satisfying \textup{(L1)} and \textup{(L2)} of Theorem~\ref{prop:algL-mub}.
\end{theorem}
\begin{proof}
Suppose that $D$ is a Lawson compact algebraic $L$-domain, and let
$X$ be its canonical IG-closure space. By
Lemma~\ref{lem:wang}(4), the map $\eta_D\colon D\to\G(X)$,
$d\mapsto\cl d$, is an order isomorphism, and $X$ is locally
consistent. Since $D$ is algebraic, so is $\G(X)$, and
Lemma~\ref{lem:compacts} shows that $X$ is dense. As in the proof
of Theorem~\ref{thm:repr_domain}, the order isomorphism $\eta_D$
is a homeomorphism for the Lawson topologies. Hence $\G(X)$ is
Lawson compact, and Theorem~\ref{prop:algL-mub} shows that
$X$ satisfies \textup{(L1)} and \textup{(L2)}.

Conversely, suppose that $D$ is order-isomorphic to $\G(X)$ for a
dense, locally consistent IG-closure space $X$ satisfying
\textup{(L1)} and \textup{(L2)} of
Theorem~\ref{prop:algL-mub}. By Lemma~\ref{lem:compacts},
$\G(X)$ is algebraic, and by Lemma~\ref{lem:wang}(3) it is an
$L$-domain. Theorem~\ref{prop:algL-mub} then shows that
$\G(X)$ is Lawson compact. Since $D$ is order-isomorphic to
$\G(X)$, it is a Lawson compact algebraic $L$-domain.
\end{proof}

\begin{remark}\label{rem:pointed-algL}
If a domain has a least element, then that element is compact and is its unique minimal element. Thus, for a dense, locally consistent IG-closure space $X$ such that $\G(X)$ has a least element, the set $\mathcal K_X=\K(\G(X))$ has a least element, which is its unique minimal element, and \textup{(L1)} holds automatically. Together with Theorem~\ref{thm:algL-mub-repr}, this shows that a poset $D$ is a pointed Lawson compact algebraic $L$-domain if and only if it is order-isomorphic to $\G(X)$ for a dense, locally consistent IG-closure space $X$ such that $\G(X)$ has a least element and \textup{(L2)} holds.
\end{remark}

By Example~\ref{ex:CMfails}, the canonical IG-closure space of
a Lawson compact algebraic domain need not satisfy \textup{(CM)}.
We now construct a representing space from its compact elements
that does satisfy \textup{(CM)}.

\begin{theorem}\label{thm:algCM}
Let $D$ be an algebraic domain. Equip $X_D=\K(D)$ with the operation
\[
\cl A=\da_D A\cap\K(D)
\qquad(A\subseteq X_D).
\]
Then the following statements hold.
\begin{enumerate}
\item $X_D$ is a dense IG-closure space, and $D$ is order-isomorphic to $\G(X_D)$.
\item $X_D$ satisfies \textup{(CM)} if and only if $D$ is Lawson compact.
\end{enumerate}
\end{theorem}

\begin{proof}
\textup{(1)} The operation $\cl\cdot$ is monotone and idempotent, and $x\in\cl x$ for every $x\in X_D$. Let $x\in X_D$ and $M\fin\cl x$. Taking $y=x$, we have $y\in\cl x$ and $M\subseteq\cl y$. Thus \textup{(In)} holds, and $X_D$ is an IG-closure space. Moreover, if $b\in\cl a$, then $a\in\cl a$ and $\{b,a\}\subseteq\cl a$. Hence $X_D$ is dense.

For each $d\in D$, put
\[
U_d=\da_D d\cap\K(D).
\]
Since $D$ is algebraic, $U_d$ is nonempty and directed. For every $M\fin U_d$, there exists $k\in U_d$ such that $m\leq k$ for every $m\in M$. Consequently,
\[
M\subseteq\cl k=\da_D k\cap\K(D)\subseteq U_d.
\]
Thus $U_d\in\G(X_D)$.

Now let $U\in\G(X_D)$. By Lemma~\ref{lem:gd}(1), $U$ is nonempty. For any $x,y\in U$, Definition~\ref{def:gd} gives $z\in U$ such that $\{x,y\}\subseteq\cl z$. This shows that $x,y\leq z$, so $U$ is directed.

Moreover, $U$ is a lower set of $\K(D)$. Indeed, if $k\in\K(D)$ and $k\leq x$ for some $x\in U$, then
\[
k\in\cl x\subseteq U,
\]
where the inclusion follows from Lemma~\ref{lem:gd}(2).

Put $d=\bigvee U$. We have $U\subseteq U_d$. Conversely, let $k\in U_d$. Then $k\in\K(D)$ and $k\leq d=\bigvee U$. Since $k$ is compact and $U$ is directed, there exists $u\in U$ such that $k\leq u$. As $U$ is a lower set of $\K(D)$, we have $k\in U$. Thus
\[
U=U_d=\da_D d\cap\K(D).
\]

The map $d\mapsto U_d$ is therefore surjective from $D$ onto $\G(X_D)$. It is also order-preserving. If $U_d\subseteq U_e$, then algebraicity gives
\[
d=\bigvee U_d\leq\bigvee U_e=e.
\]
Hence $d\leq e$ if and only if $U_d\subseteq U_e$, and the map
\[
D\longrightarrow\G(X_D),\qquad d\longmapsto\da_D d\cap\K(D),
\]
is an order isomorphism.

\textup{(2)} For $k,l\in\K(D)$, we have $\cl k=\da_D k\cap\K(D)$ and $\cl l=\da_D l\cap\K(D)$. Hence $k\leq l$ implies $\cl k\subseteq\cl l$. Conversely, if $\cl k\subseteq\cl l$, then $k\in\cl k\subseteq\cl l$, so $k\leq l$.

Since $\mathcal B_{X_D}=\set{\cl k}{k\in\K(D)}$, the map $k\mapsto\cl k$ is an order isomorphism from $\K(D)$ onto $\mathcal B_{X_D}$. By Proposition~\ref{prop:cm} and this isomorphism, $X_D$ satisfies \textup{(CM)} if and only if $\K(D)$ has property $M$. By Lemma~\ref{lem:jung-algebraic}, this is equivalent to $D$ being Lawson compact.
\end{proof}

\begin{theorem}\label{thm:algCM-repr}
A poset $D$ is a Lawson compact algebraic domain if and only if it is order-isomorphic to $\G(X)$ for a dense IG-closure space $X$ satisfying \textup{(CM)}.
\end{theorem}

\begin{proof}
Suppose that $D$ is a Lawson compact algebraic domain. By Theorem~\ref{thm:algCM}, the space $X_D$ is dense, satisfies \textup{(CM)}, and represents $D$.

Conversely, suppose that $D$ is order-isomorphic to $\G(X)$ for a dense IG-closure space $X$ satisfying \textup{(CM)}. Then $\G(X)$ is algebraic by Lemma~\ref{lem:compacts} and Lawson compact by Corollary~\ref{cor:CM}. Hence $D$ is a Lawson compact algebraic domain.
\end{proof}

For an algebraic domain $D$, Lawson compactness is equivalent
to the existence of a representation of $D$ in each of the
following three forms. The first uses a dense IG-closure space
satisfying \textup{(A1)} and \textup{(A2)}. The second uses an
IG-closure space equipped with a subset satisfying \textup{(B)}
and \textup{(M$_S$)}. The third uses a dense IG-closure space
satisfying \textup{(CM)}.

\section{Categorical equivalences}\label{sec:cat}
The categorical equivalences in this section are based on the
correspondence between approximable mappings of IG-closure spaces
and Scott-continuous maps of domains. We begin by recalling the
definition of an approximable mapping.

\begin{definition}[{\cite[Definition~6.1]{wangli}}]\label{def:approximable}
Let $X$ and $Y$ be IG-closure spaces with operations
$\cl{\cdot}_X$ and $\cl{\cdot}_Y$, respectively. A relation
$\Theta\subseteq X\times Y$ is called an \emph{approximable mapping}
if it satisfies the following conditions.
\condition{(AM1)}{For every $x\in X$, there exists $y\in Y$ such that $x\Theta y$.}
\condition{(AM2)}{For all $x\in X$ and $y,y'\in Y$, if $x\Theta y$ and $y'\in\cl y_Y$, then $x\Theta y'$.}
\condition{(AM3)}{For all $x,z\in X$ and $y\in Y$, if $x\in\cl z_X$ and $x\Theta y$, then $z\Theta y$.}
\condition{(AM4)}{For every $x\in X$ and every $M\fin Y$, if $x\Theta m$ for all $m\in M$, then there exist $z\in\cl x_X$ and $y'\in Y$ such that $z\Theta y'$ and $M\subseteq\cl{y'}_Y$.}
\end{definition}

The identity approximable mapping on an IG-closure space $X$ is
given by
\[
x\mathrel{\mathrm{id}_X}y
\quad\Longleftrightarrow\quad
y\in\cl x_X
\qquad(x,y\in X).
\]
Composition is ordinary relational composition. With these
identities and composition, IG-closure spaces and approximable
mappings form a category \cite[Proposition~6.6]{wangli}, which we
denote by $\mathbf{IG}$. Let $\mathbf{DOM}$ denote the category of
domains and Scott-continuous maps.

\begin{lemma}[{\cite[Theorems~6.5 and~6.9]{wangli}}]\label{lem:catknown}
The functor $\Fun\colon\mathbf{IG}\to\mathbf{DOM}$ given by
\[
\Fun(X)=\G(X),
\qquad
\Fun(\Theta)=\phi_\Theta
\]
is an equivalence of categories. For an approximable mapping
$\Theta\colon X\to Y$, $\phi_\Theta\colon\G(X)\to\G(Y)$ is
given by
\begin{equation}\label{eq:phi-theta}
\phi_\Theta(U)=\set{y\in Y}{x\Theta y\text{ for some }x\in U}
\qquad(U\in\G(X)).
\end{equation}
Conversely, a Scott-continuous map $f\colon\G(X)\to\G(Y)$ determines
an approximable mapping $\Theta_f\colon X\to Y$ given by
\begin{equation}\label{eq:theta-f}
x\Theta_f y\quad\Longleftrightarrow\quad y\in f(\cl x_X)
\qquad(x\in X,\ y\in Y).
\end{equation}
For each pair of IG-closure spaces $X$ and $Y$, the assignments
$\Theta\mapsto\phi_\Theta$ and $f\mapsto\Theta_f$ are mutually
inverse.
\end{lemma}

We consider the following full subcategories of $\mathbf{IG}$.
\begin{itemize}
\item $\mathbf{LCIG}$ consists of the IG-closure spaces satisfying
\textup{(C1)} and \textup{(C2)}.
\item $\mathbf{CMIG}$ consists of the IG-closure spaces $X$
satisfying \textup{(CM)} and $\G(X)=\mathcal B_X$.
\item $\mathbf{AIG}$ consists of the dense IG-closure spaces
satisfying \textup{(A1)} and \textup{(A2)}.
\item $\mathbf{CIG}$ consists of the dense IG-closure spaces
satisfying \textup{(CM)}.
\item $\mathbf{LCAIG}$ consists of the dense, locally consistent
IG-closure spaces satisfying \textup{(L1)} and \textup{(L2)}.
\end{itemize}

We also consider the following full subcategories of
$\mathbf{DOM}$.
\begin{itemize}
\item $\mathbf{LCD}$ consists of the Lawson compact domains.
\item $\mathbf{MD}$ consists of the domains $D$ with property $M$.
\item $\mathbf{LCAD}$ consists of the Lawson compact algebraic
domains.
\item $\mathbf{LCALG}$ consists of the Lawson compact algebraic
$L$-domains.
\end{itemize}

\begin{theorem}\label{thm:catdomain}
The functor $\Fun$ restricts to the following equivalences:
\begin{enumerate}
\item $\mathbf{LCIG}\simeq\mathbf{LCD}$;
\item $\mathbf{CMIG}\simeq\mathbf{MD}$;
\item $\mathbf{AIG}\simeq\mathbf{LCAD}$ and $\mathbf{CIG}\simeq\mathbf{LCAD}$;
\item $\mathbf{LCAIG}\simeq\mathbf{LCALG}$.
\end{enumerate}
\end{theorem}

\begin{proof}
	In each part, the category of spaces and the category of domains
	are full subcategories of $\mathbf{IG}$ and $\mathbf{DOM}$,
	respectively. Thus, whenever $\Fun$ maps the former into the
	latter, its restriction is full and faithful by
	Lemma~\ref{lem:catknown}.
	
	For \textup{(1)}, Theorem~\ref{thm:main} shows that $\Fun$
	maps $\mathbf{LCIG}$ into $\mathbf{LCD}$. The restriction is
	essentially surjective by Theorem~\ref{thm:repr_domain}.
	
	For \textup{(2)}, Proposition~\ref{prop:cm} shows that $\Fun$
	maps $\mathbf{CMIG}$ into $\mathbf{MD}$. The restriction is
	essentially surjective by Theorem~\ref{thm:whole} and
	Lemma~\ref{lem:wang}(4).
	
	For \textup{(3)}, Lemma~\ref{lem:compacts},
	Theorem~\ref{thm:algebraic} and Corollary~\ref{cor:CM} show
	that $\Fun$ maps both $\mathbf{AIG}$ and $\mathbf{CIG}$ into
	$\mathbf{LCAD}$. The respective restrictions are essentially
	surjective by Theorem~\ref{thm:alg_domain} and
	Theorem~\ref{thm:algCM}.
	
	For \textup{(4)}, Lemma~\ref{lem:compacts},
	Lemma~\ref{lem:wang}(3) and Theorem~\ref{prop:algL-mub} show that
	$\Fun$ maps $\mathbf{LCAIG}$ into $\mathbf{LCALG}$. The
	restriction is essentially surjective by
	Theorem~\ref{thm:algL-mub-repr}.
	
	Therefore all the stated restrictions are equivalences.
\end{proof}

Finally, we extend the representation in Theorem~\ref{thm:selected} to an equivalence of categories. 
\begin{itemize}
	\item
Let $\mathbf{SIG}$ be the category whose objects are triples
$(X,\cl{\cdot},S)$, where $(X,\cl{\cdot})$ is an IG-closure space and
$S\subseteq X$ satisfies \textup{(B)} and \textup{(M$_S$)}. The
morphisms are approximable mappings between the underlying
IG-closure spaces, with no requirement on
the selected subsets. Identities and composition are those of
$\mathbf{IG}$.
\item Let $\mathbf{MB}$ be the category whose objects are triples
$(D,\leq,\mathcal B)$, where $(D,\leq)$ is a domain and $\mathcal B$ is a basis of $D$
with property $M$. The morphisms are Scott-continuous maps between the
underlying domains. Identities and composition are those of
$\mathbf{DOM}$.
\end{itemize}

\begin{theorem}\label{thm:selectedcat}
The categories $\mathbf{SIG}$ and $\mathbf{MB}$ are equivalent.
\end{theorem}

\begin{proof}
For each object $(X,\cl{\cdot},S)$ and each morphism
$\Theta$ of $\mathbf{SIG}$, put
\[
\Funs(X,\cl{\cdot},S)
=
\bigl(\G(X),\subseteq,\mathcal B_S\bigr),
\qquad
\Funs(\Theta)=\phi_\Theta,
\]
where $\phi_\Theta$ is the Scott-continuous map determined by
$\Theta$ in Lemma~\ref{lem:catknown}.
By Lemmas~\ref{lem:subbasis} and~\ref{lem:ms}, $\mathcal B_S$ is a
basis of $\G(X)$ with property $M$.
Hence $\bigl(\G(X),\subseteq,\mathcal B_S\bigr)$ is an object of
$\mathbf{MB}$. By Lemma~\ref{lem:catknown}, $\phi_\Theta$ is Scott
continuous, and $\Funs$ preserves identities and composition. Thus
$\Funs$ is a functor.

Let $(X,\cl{\cdot},S)$ and $(Y,\cl{\cdot},T)$ be objects of
$\mathbf{SIG}$. To show that $\Funs$ is faithful, suppose that
$\Theta_1,\Theta_2\colon X\to Y$ are approximable mappings such that
$\Funs(\Theta_1)=\Funs(\Theta_2)$. Then
$\phi_{\Theta_1}=\phi_{\Theta_2}$, so $\Theta_1=\Theta_2$ by
Lemma~\ref{lem:catknown}.

To show that $\Funs$ is full, let
\[
f\colon\bigl(\G(X),\subseteq,\mathcal B_S\bigr)
\longrightarrow
\bigl(\G(Y),\subseteq,\mathcal B_T\bigr)
\]
be a morphism in $\mathbf{MB}$. Then $f\colon\G(X)\to\G(Y)$ is Scott
continuous. By Lemma~\ref{lem:catknown}, there exists an
approximable mapping $\Theta_f\colon X\to Y$ such that
$\phi_{\Theta_f}=f$. By the definition of morphisms in
$\mathbf{SIG}$, $\Theta_f$ is a morphism from $(X,\cl{\cdot},S)$ to
$(Y,\cl{\cdot},T)$. Hence $\Funs(\Theta_f)=f$, and $\Funs$ is full.

Finally, let $(D,\leq,\mathcal B)$ be an object of $\mathbf{MB}$. Let $X$ be
the canonical IG-closure space of $D$ and put $S=\mathcal B$. By
Lemma~\ref{lem:wang}(4), the map
\[
\eta_D\colon D\longrightarrow\G(X),
\qquad
d\longmapsto\cl d,
\]
is an order isomorphism, and $\eta_D(\mathcal B)=\mathcal B_S$. Therefore
$\mathcal B_S$ is a basis of $\G(X)$ with property $M$. By Lemmas~\ref{lem:subbasis} and~\ref{lem:ms},
the subset $S$ satisfies \textup{(B)} and \textup{(M$_S$)}. Thus
$(X,\cl{\cdot},S)$ is an object of $\mathbf{SIG}$.

Since $\eta_D$ and its inverse are Scott continuous, $\eta_D$ is an
isomorphism in $\mathbf{MB}$ from $(D,\leq,\mathcal B)$ to
$\bigl(\G(X),\subseteq,\mathcal B_S\bigr)$. Hence $\Funs$ is
essentially surjective. Therefore $\Funs$ is an equivalence of
categories.
\end{proof}
\section*{Statements and Declarations}
\noindent\textbf{Competing Interests.} The author declares that there
are no competing interests.

\section*{Acknowledgements}
{\emergencystretch=2em
This work is supported by the National Natural Science Foundation of China (Grant No.~12571489).\par}

\end{document}